\documentclass[10pt]{amsart}

\usepackage[matrix,arrow,curve,frame]{xy}
\usepackage{amsmath,amsthm,amssymb,enumerate}
\usepackage{latexsym}
\usepackage{amscd}

\newtheorem{thm}[subsection]{Theorem}
\newtheorem{defn}[subsection]{Definition}
\newtheorem{prop}[subsection]{Proposition}

\newtheorem{cor}[subsection]{Corollary}
\newtheorem{lemma}[subsection]{Lemma}

\newtheorem{remark}[subsection]{Remark}

\theoremstyle{definition}

\newcommand{\cat}{\mathcal}
\newcommand{\lra}{\longrightarrow}

\newcommand{\R}{\mathbb R}
\newcommand{\Q}{\mathbb Q}
\newcommand{\Z}{\mathbb Z}
\newcommand{\N}{\mathbb N}
\newcommand{\C}{\mathbb C}
\newcommand{\CP}{\mathbb P}
\newcommand{\T}{\mathbb T}

\newcommand{\hT}{{\mathbb T}}

\newcommand{\I}{{\cat I}}
\newcommand{\Jj}{{\cat J}}
\newcommand{\Ll}{{\cat L}}
\newcommand{\Oo}{{\cat O}}

\newcommand{\KS}{{\cat KS}}
\newcommand{\Mm}{{\cat M}}
\newcommand{\A}{{\cat A}}
\newcommand{\K}{\cat K}

\newcommand{\Ga}{\Gamma}
\newcommand{\Om}{\Omega}
\newcommand{\om}{\omega}
\newcommand{\al}{\alpha}
\newcommand{\la}{\lambda}
\newcommand{\Si}{\Sigma}

\newcommand{\db}{\bar{\partial}}

\DeclareMathOperator{\im}{im}
\DeclareMathOperator{\Aut}{Aut}

\DeclareMathOperator{\hocolim}{hocolim}
\DeclareMathOperator{\holim}{holim}

\DeclareMathOperator{\Map}{Map}

\DeclareMathOperator{\Hom}{Hom}

\DeclareMathOperator{\Diff}{Diff}
\DeclareMathOperator{\FDiff}{FDiff}

\DeclareMathOperator{\Symp}{Symp}
\DeclareMathOperator{\Lie}{Lie}

\DeclareMathOperator{\Hol}{Hol}
\DeclareMathOperator{\Iso}{Iso}
\DeclareMathOperator{\Det}{Det}
\DeclareMathOperator{\Sym}{Sym}

\begin{document}

\title[Compatible complex structures]
{Compatible complex structures on symplectic
rational ruled surfaces}
\author{Miguel Abreu}
\address{Departamento de Matem\'atica, Instituto Superior T\'ecnico,   
Av. Rovisco Pais,\newline\indent 1049-001 Lisboa, Portugal}
\email{mabreu@math.ist.utl.pt, ggranja@math.ist.utl.pt}
\author{Gustavo Granja}
\author{Nitu Kitchloo}
\address{Department of Mathematics, University of California, San
  Diego, USA}
\email{nitu@math.ucsd.edu}
\thanks{Nitu Kitchloo is supported in part by NSF through grant DMS
  0436600. Miguel Abreu and Gustavo Granja are supported in part by 
  FCT through program POCTI-Research Units Pluriannual Funding Program 
  and grants POCTI/MAT/57888/2004 and POCTI/MAT/58497/2004.}

\date{\today}


\dedicatory{In fond memory of Raoul Bott}

{\abstract

\noindent
In this paper we study the topology of the space $\I_\omega$ of 
complex structures compatible with a fixed symplectic form $\omega$, 
using the framework of Donaldson. By comparing our analysis of the
space $\I_\omega$ with results of McDuff on the space $\cat J_\omega$ 
of compatible almost complex structures on rational ruled surfaces, 
we find that $\I_\omega$ is contractible in this case.

We then apply this result to study the topology of the
symplectomorphism group of a rational ruled surface, extending results
of Abreu and McDuff.}

\maketitle

\section{Introduction}

\noindent
The work of Gromov on $J$-holomorphic curves \cite{Gr} has provided tools for 
understanding the topology of the group of symplectomorphisms of
certain four-dimensional symplectic manifolds. Gromov began the study of 
the group of symplectomorphisms on $S^2 \times S^2$ with the standard 
symplectic form $\sigma \oplus \sigma$. He showed that this group is 
homotopy equivalent to a semidirect product of $\Z/2$ with $SO(3)
\times SO(3)$. Later Abreu \cite{Ab} continued this study by analyzing 
the structure of the group of symplectomorphism on $S^2 \times S^2$
with symplectic form $\omega = \lambda \sigma \oplus \sigma, \; 1 <
\lambda \leq 2$. This work was extended by Abreu-McDuff \cite{AM} to 
arbitrary $\lambda$, leading to a complete calculation of the rational 
cohomology ring of the classifying space of the symplectomorphism group 
(modulo a mistake which is corrected below in Theorem \ref{Main2}). There are
corresponding results for the structure of the symplectomorphism group
of the non trivial $S^2$-bundle over $S^2$ (for an arbitrary
symplectic form). The basic idea in the work mentioned above is to
analyze the action of the symplectomorphism group on the contractible
space of compatible almost complex structures $\Jj_{\omega}$. 

In \cite{Do} Donaldson showed that the action of the symplectomorphism
group on $\cat J_{\omega}$ is Hamiltonian, and the moment map for this 
action is the Hermitian scalar curvature of the corresponding almost 
K\"ahler metric. He also showed that, restricted to the space
$\I_{\omega}$ of compatible integrable structures, this action fits
into the general framework of infinite dimensional geometric invariant 
theory, going back to Atiyah and Bott \cite{AB1} (see \cite{AK1} for more
discussion of this point of view). Therefore, in principle, the norm square of the moment
map should induce a stratification of $\I_{\omega}$ with critical 
points being the extremal K\"ahler metrics compatible with the 
symplectic form. By work of Calabi \cite{C1}, for rational ruled surfaces these metrics 
correspond to a finite collection of 
Hirzebruch surfaces. In particular, this suggests 
that each stratum in $\I_{\omega}$ should be homotopy equivalent to a 
single orbit of the symplectomorphism group. A stratification of 
$\cat J_{\omega}$ with similar properties had been established by 
Abreu and studied in detail by McDuff \cite{McD1}. This indicated that 
the hypothetical stratification on $\I_{\omega}$ determined by the 
moment map could be the one induced via the inclusion 
$\I_{\omega} \subset \Jj_{\omega}$. 

We begin our study of the space $\I_{\omega}$ by analyzing the 
stratification induced by the above inclusion. We show that each 
stratum $V$ contains an orbit of the symplectomorphism group, 
corresponding to a Hirzebruch surface $F$ with a standard K\"ahler  
metric, which is weakly equivalent to it.
The stratification of $\I_{\omega}$ has the advantage that its 
gluing data can be easily understood via Kodaira-Spencer deformation theory. 
By comparing our stratification of $\I_{\omega}$ with that of Abreu-McDuff, we prove:

\begin{thm} \label{Main1}
The inclusion of the space of compatible integrable complex structures  
into the space of all compatible almost complex structures, 
$\I_{\omega}(M) \subset \Jj_{\omega}(M)$, is a weak homotopy
equivalence for a rational ruled surface $M$. In particular, the space 
$\I_\omega(M)$ is weakly contractible.
\end{thm}

\noindent
As far as we are aware, and besides the obvious case of real
$2$-dimensional surfaces, this is the first example where the topology 
of the space of compatible integrable complex structures on a
symplectic manifold has been understood. 

Recall that any rational ruled surface is diffeomorphic to either $S^2
\times S^2$, the trivial $S^2$-bundle over $S^2$, or
$S^2\tilde{\times}S^2$, the non trivial $S^2$-bundle over $S^2$. Work
of Taubes, Liu-Li and Lalonde-McDuff (see~\cite{LM} for detailed
references) implies that any symplectic form on one of these smooth
manifolds is ``standard'', i.e. diffeomorphic to a scalar
multiple of $\om_\la =\la\sigma\oplus\sigma$, $1\leq\la\in\R$, on $S^2
\times S^2$, or to any chosen symplectic form $\om_\la$, $0<\la\in\R$,
on $S^2\tilde{\times}S^2$ such that $[\om_\la](E) = \la$ and
$[\om_\la](F) = 1 $, where $E$ denotes the homology class of the
exceptional divisor under the natural identification
$S^2\tilde{\times}S^2 \cong \CP^2 \sharp \overline{\CP^2}$ and $F$
denotes the homology class of the fiber.

Complex deformation theory gives us a good understanding of the 
way the strata of $\I_\omega$ glue and allows us to express the symplectomorphism
group as an iterated homotopy pushout of certain compact subgroups
(see Theorem \ref{maintopological}). This, in turn immediately gives 
the integral cohomology groups of the classifying space of the symplectomorphism groups of
rational ruled surfaces. Let $G_\lambda$ denote the group of 
symplectomorphisms of $S^2 \times S^2$ with  symplectic form 
$\om_\la = \lambda \sigma \oplus \sigma$, where $1<\lambda\in\R$ lies 
between the integers $0 < \ell < \lambda \leq \ell + 1$. Let $BG_\lambda$ denote the
classifying space of $G_\la$. We have:

\begin{thm} \label{cohgps}
The integral cohomology groups of $BG_\lambda$ in the untwisted case are given by:
\[
H^*(BG_\lambda; \Z) =  H^*(BSO(3)\times BSO(3));\Z) \oplus \bigoplus_{i=1}^\ell \Sigma^{4i-2} 
H^*(BS^1\times BSO(3);\Z),
\]
where $\Sigma$ denotes the suspension of graded abelian groups.
\end{thm}

\noindent
The relevance of the compact Lie groups appearing in the previous statement
 will become clear in sections 3--5. If we work away from the prime 2, we can compute 
 the ring structure:

\begin{thm} \label{Main2}
The cohomology of $BG_\lambda$ in the untwisted case with coefficients in the ring $R= \Z[1/2]$ is 
given by the following free module over the ring $R[x,y]$ on generators 
$b_i$, $a_j$, $ 0 \leq i < \ell$, $0 \leq j \leq \ell$:
\[ H^*(BG_\lambda;R) = R[x,y]\langle a_0,b_0,a_1,b_1,a_2,\ldots a_\ell \rangle,  \]
where $a_0=1$, the degree of the elements $x,y$ is $4$, the degree of $b_k$ is $4k+2$, 
and that of $a_k$ is $4k$. Moreover, as a ring, we may identify 
$H^*(BG_\lambda,R)$ as the subring of 
\[ 
H^*(BG_\lambda;\Q) = \frac{\Q[x,y,z]}{ \langle z\prod_{i=1}^\ell
(z^2+i^4 x - i^2 y) \rangle} 
\]
where the degree of $z$ is $2$, and the elements $a_k$ and $b_k$ are identified 
respectively with the elements
\[ 
\frac{z^2}{(2k)!}\prod_{i=1}^{k-1}(z^2+ i^4 x-i^2 y),\quad \text{and } \quad\frac{z}{(2k+1)!}\prod_{i=1}^k(z^2+ i^4 x-i^2 y).
\]

\end{thm}
\noindent
If $G_\lambda$ denotes the group of symplectomorphisms of 
$S^2 \tilde{\times} S^2$, with symplectic form $\om_\la$ as above,
where $0<\lambda\in\R$ lies between the integers 
$0 \leq \ell < \lambda \leq \ell + 1$, we have:

\begin{thm} \label{cohgps2}
The cohomology groups of the space $BG_{\lambda}$ in the twisted case are given by:
\[
H^*(BG_\lambda; \Z) =  \bigoplus_{i=0}^\ell \Sigma^{4i} 
H^*(BU(2);\Z).
\]
\end{thm}

\noindent
In particular we see that the cohomology of $BG_\lambda$ is torsion free in the twisted case.

\begin{thm} \label{Main3}
The rational cohomology ring of $BG_{\lambda}$ in the twisted case is given by:
\[ 
H^*(BG_{\lambda};\Q) = \frac{\Q[x,y,z]}{ \langle \prod_{i=0}^{\ell}(- z^2+(2i+1)^4 x - (2i+1)^2 y) \rangle} 
\]
where the degree of the class $z$ is 2, and that of the classes $x$ and $y$ is 4.

\end{thm}

\bigskip
\noindent
We bring to the attention of the reader the relations in the rational
cohomology of $BG_\lambda$. In \cite{AM}, both in the twisted and untwisted cases, the decomposable term $z^2$ was missing in each of the factors that make up the relation. In addition, in the twisted case, the nondecomposable terms in the factors have powers of odd integers as 
coefficients instead of all integers. The source of these inaccuracies in \cite{AM} is discussed in Remark \ref{Dusamistake}.

Regarding previous papers on the topology of symplectomorphism groups of rational ruled 
surfaces, the results of this paper depend logically only on the analysis of the 
stratification of $\Jj_\omega$ in \cite{McD1}, the computation of the homotopy type
of the strata in \cite{Ab} (or \cite{AM}) and (just for the computation of the ring structures) 
the fact,  proved in \cite{McD2}, that the homotopy colimit of 
certain inclusions $BG_\lambda \to BG_{\lambda+\epsilon}$ is the classifying space of the group of fiberwise diffeomorphisms.

\subsection{Organization of the Paper:}

\noindent

\noindent
The body of the paper is divided into four sections and three appendices. 

Section 2 is further divided into four subsections. In the first part, 
we begin by developing our framework to study complex deformation
theory of a complex 4-manifold $M$. In the next part, we specialize to
a symplectic 4-manifold $(M,\omega)$. Let $\Jj_{\omega}$ be the space
of almost complex structures on $M$ compatible with $\omega$, and let 
$\I_{\omega}$ be the compatible integrable structures. Given a
K\"ahler structure $J_0 \in \I_{\omega}$, we interpret the normal
bundle in $\I_{\omega}$ to the space of equivalent K\"ahler
structures, in terms of deformation theory. In the next subsection, 
we study the action of the diffeomorphism group on the space of 
compatible complex structures for an arbitrary symplectic manifold 
$(M,\omega)$. Using this we show that, under certain conditions,
the space of diffeomorphic compatible structures is homotopy equivalent
 to an orbit of the symplectomorphism group. Finally, in the last part
of this section, we analyze the inclusion $\iota : \I_{\omega} \subset \Jj_{\omega}$. 
The space $\Jj_{\omega}$ admits a map from the moduli space $\cat M$ 
of $J$-holomorphic curves. We derive sufficient cohomological
conditions for $\cat M$ to meet $\iota $ transversally. In the case of 
rational ruled surfaces, this implies that the map $\iota$ is
transverse to the strata introduced by Abreu. Some technical 
lemmas in section 2 have been banished to Appendices \ref{app:A} 
and \ref{app:C}. We feel that many of the results in this section should 
be well known to experts but we have not been able to find convenient 
references in the literature.

In section 3, we study rational ruled surfaces in detail. This section
is divided into two parts. In the first part, we describe a
stratification on the space  $\Jj_{\omega}$ of compatible almost
complex structures, for a fixed symplectic form $\omega$, previously
studied in~\cite{Ab,AM, McD1}. In the second part, we use results from
the previous section to show that this induces an analogous
stratification on the corresponding space $\I_{\omega}$ of the
compatible complex structures. We then conclude with the proof of
Theorem~\ref{Main1}.  

Section 4 is dedicated to studying the deformation theory of
Hirzebruch surfaces. We review the construction of Hirzebruch 
surfaces by K\"ahler reduction and use it to identify their K\"ahler 
automorphism groups. We also relate standard bases for the maximal tori of the automorphism groups
to the bases yielding standard Delzant polygons as images of the moment map. We then 
use the fixed point formula for elliptic complexes of Atiyah and Bott to
determine the isotropy representations of the symplectomorphism group 
on the normal bundle to the various strata in $\I_{\omega}$. 
These results are applied in the next section to study the topology of the
symplectomorphism group.

In section 5, we use the results of sections 3 and 4 to express the 
classifying space of the symplectomorphism group as a finite iterated
homotopy pushout of classifying spaces of certain K\"ahler isometry
subgroups (leaving some technicalities for Appendix \ref{quasi}). 
Theorems \ref{cohgps} and \ref{cohgps2} follow immediately.
We then use the loop maps from $G_\lambda$ to the group of fiberwise
diffeomorphisms defined in \cite{McD2}
together with the classification of Hamiltonian $S^1$-actions on 
four-manifolds \cite{Ka} to compute the rational cohomology rings described in 
Theorems \ref{Main2} and \ref{Main3} (see Theorem \ref{Main23} and Remark 
\ref{concrem}). A key step is understanding the relation between the cohomology
of the various K\"ahler isometry groups (Proposition \ref{mapsoncohomology}). 
Finally, we use the computation 
in \cite{HHH} of the $T^2$-equivariant cohomology groups of $\Omega SU(2)$ 
to find the cohomology ring of the classifying space of the fiberwise
diffeomorphism group away from $2$ and use this to complete the proof of
Theorem \ref{Main2}.

\subsection{Acknowledgments:}

\noindent

\noindent
The authors would like to thank Dusa McDuff for various helpful conversations and a referee for carefully  reading the paper and suggesting improvements to the exposition. The first author would like to thank
Vestislav Apostolov for a useful reference regarding the material in Appendix \ref{app:A}. The second author would like to thank Sue Tolman for a very helpful discussion. The third author would like to thank the Department of Mathematics of Instituto Superior T\'ecnico for its hospitality while this work was being conducted. 

\subsection{Conventions:}

\noindent

\noindent
Throughout this paper, we work with infinite dimensional Fr\'echet
manifolds which are locally modeled on an inverse limit
of Banach spaces. Standard theorems such as the inverse function
theorem do not hold automatically in the Fr\'echet setting. 
Therefore, we need to say a few words about the context in which the 
transversality arguments in this paper need to be interpreted. 
All the Fr\'echet manifolds we work with can naturally be interpreted 
as inverse limits of Banach manifolds. For example, the space of smooth 
sections of a bundle over a smooth manifold is the intersection over 
$k$ of the corresponding Banach manifolds of $C^k$-sections. 
For each individual $k$, all the transversality arguments we use hold 
in the infinite dimensional context. The validity of the results stated 
in the smooth setting should therefore be interpreted as the validity of 
the corresponding result for each Banach manifold indexed by $k$. 
Statements about the homotopy type of the corresponding Fr\'echet manifold 
can be derived from the fact that the successive inclusions between the  
Banach manifolds are weak equivalences (see \cite{P}). 
We shall illustrate this with an example in Remark~\ref{rmk:cdef}.

\section{General facts on compatible complex structures}

\noindent
The goal of this section is to set up a geometric framework and
establish some facts regarding the space of compatible integrable 
complex structures $\I_\om$ and its inclusion
in the space of compatible almost complex structures $\Jj_\om$ 
on a symplectic $4$-manifold $(M,\om)$. 

We begin by describing cohomological conditions under which the space $\I_\om$ is a submanifold of
 $\Jj_\om$ (see Theorem \ref{thm:intsubacs}). We then show that for 
 a given $J \in \I_\om$, the intersection $(\Diff(M)\cdot J) \cap \I_\om$ is weakly equivalent to 
 the orbit of $J$ under the symplectomorphism group, as long as the K\"ahler isometry group
 of $(J,\omega)$ is a deformation retract of the complex automorphisms of $J$ preserving
 the cohomology class of $\om$ (see Corollary \ref{cor:intorbit}). Finally, we find conditions under which the projection map from the space of $J$-holomorphic curves 
to $\Jj_\om$ is transverse to $\I_\om$. In the case of rational ruled surfaces this will imply that Abreu and McDuff's stratification of $\Jj_\omega$ induces a stratification of $\I_\om$ (see Theorem \ref{thm:transversal}).

\subsection{Complex structures:}

\noindent

\noindent 
In this subsection we review the classical deformation
theory of Kodaira-Spencer (see \cite{Ko}) from the point of view that we
will adopt in the next subsection to study deformations of compatible
complex structures.

\noindent
Given an almost complex manifold $(M,J)$, the Nijenhuis tensor 
$N_J \in \Om_J^{0,2} (TM) = \Om_J^{0,2} (M)\otimes \Om^{0}(TM)$
is a $(0,2)$-form on $(M,J)$ with values in $(TM,J)$ that measures the
non-integrability of $J$ (see Definition~\ref{defn:nijenhuis} in 
Appendix~\ref{app:A}). If one considers the space $\Jj$ of all almost
complex structures on the manifold $M$, the Nijenhuis tensor $N$ can
be seen as a section of the natural vector bundle $\Om^{0,2}(TM)$ over
$\Jj$, whose fiber over a point $J\in\Jj$ is $\Om^{0,2}_J (TM)$:
\[ N:\Jj  \to\Om^{0,2}(TM)\;, \quad J  \mapsto (J, N_J). \]
The space $\I$ of integrable complex structures on $M$ is
the zero-set of this Nijenhuis section $N$. As usual, it will be a
submanifold of $\Jj$ if the Nijenhuis section $N$ is transversal to
the zero section. 

The vector bundle $\Om^{0,2}(TM)$ is a canonical summand of the
trivial bundle over $\Jj$ with fiber $\Om^2(TM)$:
\[
\Om^{0,2}(TM) \oplus \left(\Om^{2,0}(TM)\oplus\Om^{1,1}(TM) \right) 
\equiv \Om^2 (TM) \times \Jj \to \Jj\,.
\]
This means in particular that $\Om^{0,2}(TM)$ has a natural connection
$\nabla$, given by projection of the trivial connection on $\Om^2(TM)
\times \Jj$:
\[
\nabla\cdot = (d\cdot)^{0,2}\,.
\]
Since $T\Jj\cong \Om^{0,1}(TM)$, we have that $\nabla N$ can be regarded as
a bundle map
\[
\nabla N : \Om^{0,1}(TM) \to \Om^{0,2}(TM)\,.
\]
In particular, if for a given $J\in\I$, i.e. a $J\in \Jj $ for which
$N_J\equiv 0$, the map
\[
\nabla N_J : \Om^{0,1}_J(TM) \to \Om^{0,2}_J(TM)
\]
is surjective, then the Nijenhuis section is transversal to the zero
section at this $J\in\I$.

As proved in Corollary~\ref{cor:dN02}, 
the map $\nabla N$ is essentially the $\db$-operator (see
Definition~\ref{defn:db01TM}). More precisely,
\[
\nabla N_J = (-2J)\db_J\,,\ \forall\,J\in\Jj\,.
\]
The following proposition is then immediate.
\begin{prop} \label{prop:intsubacs}
If $M$ is a $4$-dimensional manifold, $J\in\I$ is an
integrable complex structure on $M$ and the cohomology group 
$H^{0,2}_J (TM) = 0$, then $\I$ is a submanifold of $\Jj$
in the neighborhood of $J$, with tangent space
\[
T_J \I = \ker \left\{ \db:\Om^{0,1}_J(TM) \to
  \Om^{0,2}_J(TM) \right\} \subset \Om^{0,1}_J (TM) = T_J \Jj\,. 
\]
\end{prop}
\noindent
The group $\Diff(M)$ of diffeomorphisms of $M$ acts naturally on $\Jj$
via
\[
\varphi^*(J) := \left(d\varphi\right)^{-1} J
\left(d\varphi\right)\,, \ \forall\, \varphi\in\Diff(M),\, J\in\Jj\,. 
\]
The induced infinitesimal action is given by the Lie derivative
\[
\Ll:\Om^0 (TM) \to \Om^{0,1}(TM) = T\Jj
\]
which, as proved in Proposition~\ref{prop:lieJ}, can be written for a given
$X\in\Om^0 (TM)$ and at a given $J\in\Jj$ as
\[
\Ll_X J = (2J) (\db X) + \frac{1}{2} J (X \lrcorner N_J) \in
\Om^{0,1}_J (TM) = T_J \Jj\,.
\]
The action of $\Diff(M)$ on $\Jj$ preserves $\I$. Since $\db$ commutes with $J$ (see
Proposition~\ref{prop:db0TM}), the
previous formula implies that the tangent space to an
orbit $\Diff(M)\cdot J$ at a point $J\in\I$, is given by
\[
T_J (\Diff(M)\cdot J) = \im \left\{ \db:\Om^{0}_J(TM) \to
  \Om^{0,1}_J(TM) \right\}\,.
\]
It then follows from Proposition~\ref{prop:intsubacs} that the 
moduli space of infinitesimal deformations of $J\in\I$ is
given by $H^{0,1}_J (TM)$ and 
\[
T_J \I \cong T_J (\Diff(M)\cdot J) \oplus H^{0,1}_J (TM)
\,.  
\]
\begin{remark} \label{rmk:cdef}
As promised, we elaborate on the transversality argument here: 

We may see $N$ as a section of the bundle $\Omega_k^{0,2}(TM)$ 
with base being the space $\cat J_{k+1}$ (the space of $C^{k+1}$-almost 
complex structures), and fiber $\Omega_k^{0,2}(TM,J)$ (the space of 
$C^k$-forms of type (0,2) with values in $TM$) over $J \in \cat J_{k+1}$. 
This bundle supports a natural connection and the same proof shows that 
$\nabla N$ can be identified with the $\bar{\partial}$ operator
\[ 
\nabla N = \bar{\partial} : \Omega_{k+1}^{0,1}(TM,J) \lra 
\Omega_k^{0,2}(TM,J) 
\]
Let $H^{0,n}_k(TM,J)$ be the cohomology groups of the Dolbeault complex
\[ 
\Omega_{k+2}^0(TM,J) \xrightarrow{\db} \Omega_{k+1}^{0,1}(TM,J) \xrightarrow{\db}
\Omega_k^{0,2}(TM,J) \lra 0 
\]
Notice that $H_k^{0,n}(TM,J)$ is isomorphic to $H^{0,n}(TM,J)$, 
for $J \in \I$, and for all $n,k \geq 0$ since the complex above is 
a resolution of the sheaf of holomorphic vector fields
on $M$ by fine sheaves. Hence we may apply the implicit function theorem for each $k$ 
to show that the zero locus of $N$ is a Banach manifold 
$\I_{k+1} \subset \Jj_{k+1}$ modeled on the Banach space of closed 
$C^{k+1}$-forms of type $(0,1)$. 
 
Now let $\Diff_{k+2}(M)$ denote the group of $C^{k+2}$ diffeomorphisms 
of $M$. The action map of $\Diff_{k+2}(M)$ on $\Jj_{k+1}$ is $C^1$, 
so we may differentiate the action. Repeating the proof of the
previous claim, the tangent space to the orbit of the group 
$\Diff_{k+2}(M)$ at $J \in \I \cap \Jj_{k+1}$ can be identified with
the space of exact $C^{k+1}$-forms of type $(0,1)$. Therefore, the
space of infinitesimal deformations of $J$ inside $\I_{k+1}$ can be 
naturally identified with the cohomology group $H^{0,1}(TM,J)$ 
(which is independent of $k$).
\end{remark}

\subsection{Compatible complex structures:}

\noindent

\noindent
Now let $(M,\om)$ be a symplectic $4$-manifold and $\Jj_\om$ the
contractible submanifold of $\Jj$ consisting of almost complex
structures on $M$ compatible with $\om$. Given $J\in\Jj_\om$, denote
by
\[
h_J(\cdot,\cdot) \equiv \om(\cdot,J\cdot) - i \om(\cdot,\cdot)
\]
the hermitian metric on $TM$ induced by the pair $(\om,J)$. 
We may use $h_J$ to identify $T_J \Jj = \Om_J^{0,1} (TM)$ with the
space $T_J^{0,2}(M) \equiv \Om_J^{0,1}(M)\otimes \Om_J^{0,1}(M)$ of
complex $(0,2)$-tensors, via
\[
\Om_J^{0,1}(TM) \ni A(\cdot) \leftrightarrow \theta_A(\cdot,\cdot)
\equiv h_J (A\cdot, \cdot) \in T_J^{0,2}\,.
\]
Under this identification, the subspace $T_J \Jj_\om =
S\Om_J^{0,1}(TM) \subset \Om_J^{0,1}(TM) = T_J \Jj$ can be identified
with the subspace of complex symmetric $(0,2)$-tensors $S_J^{0,2}(M)
\subset T_J^{0,2}(M)$:
\[
A \in T_J \Jj_\om \Leftrightarrow AJ+JA = 0 \quad\text{and}\quad 
\om (A\cdot,\cdot) = - \om (\cdot,A\cdot) \Leftrightarrow 
\theta_A \in S_J^{0,2}(M)\,.
\]
The quotient may therefore be identified with the space of
$(0,2)$-forms on $M$:
\[
T_J \Jj / T_J \Jj_\om = \Om_J^{0,1} (TM) / S \Om_J^{0,1}(TM) \cong
T_J^{0,2}(M) / S_J^{0,2} (M) = \Om_J^{0,2} (M)\,.
\]
Given a compatible integrable complex structure $J\in\I_\om
\subset \Jj_\om$, consider the following sequence of chain complexes
with exact columns, where the above identifications are taken into
account and commutativity of the lower left corner is proved in 
Appendix~\ref{app:C}: 
\[ 
\xymatrix{  
0 \ar[r] \ar[d]  & S\Om_J^{0,1}(TM) \ar[r]^{\bar{\partial}} 
\ar[d] & \Om_J^{0,2}(TM) \ar[r] \ar[d] & 0 \\ 
\Om^0(TM) \ar[r]^{\bar{\partial}} \ar[d] & \Om_J^{0,1}(TM) 
\ar[r]^{\bar{\partial}} \ar[d] & \Om_J^{0,2}(TM) \ar[r] \ar[d] & 0 \\
\Om_J^{0,1}(M) \ar[r]^{\bar{\partial}} & \Om_J^{0,2}(M) 
\ar[r]^{\bar{\partial}} & 0 \ar[r] & 0   
} 
\]
The snake lemma yields a long exact sequence of holomorphic cohomology
groups: 
\[ 
0 \longrightarrow H^0_J(TM) \longrightarrow cl\Om_J^{0,1}(M) 
\overset{\delta}{\longrightarrow} 
clS\Om_J^{0,1}(TM) \longrightarrow  H_J^{0,1}(TM) 
\longrightarrow 
\]
\[ 
\longrightarrow H_J^{0,2}(M) \longrightarrow SH_J^{0,2}(TM)
\longrightarrow H_J^{0,2}(TM) \longrightarrow  0 
\]
where $cl\Om_J^{0,1}(M)$ and $clS\Om_J^{0,1}(TM)$ denote the closed
$(0,1)$-forms in the respective spaces and the symmetric cohomology
$SH_J^{0,2}(TM)$ is given by the cokernel of $\db$ restricted to the
symmetric $(0,2)$-tensors $S\Om_J^{0,1}(TM)$. A geometric
interpretation of the maps in the above long exact sequence can be
given as follows.

By analysing the Nijenhuis tensor and its covariant derivative as
before, one sees that $\I_\om$ is a submanifold of
$\Jj_\om$ in the neighborhood of $J\in\I_\om$ if the
symmetric cohomology group $SH_J^{0,2}(TM)$ vanishes. Under this
assumption, $clS\Om_J^{0,1}(TM)$ can be identified with the tangent
space $T_J \I_\om$, and the image of $\delta:
cl\Om_J^{0,1}(M) \to clS\Om_J^{0,1}(TM)$ identifies the tangent space
to the intersection of the orbit $(\Diff(M)\cdot J)$ with
$\I_\om$. Moreover, the image of the cokernel of $\delta$
in $H_J^{0,1}(TM)$ identifies the the moduli of infinitesimal
deformations of $J\in\I_\om \subset \I$ that can be realized in an 
$\om$-compatible way. We then have the following theorem.
\begin{thm} \label{thm:intsubacs}
If $(M,\om)$ is a symplectic $4$-dimensional manifold, 
$J\in\I_\om$ is a compatible integrable complex structure
on $(M,\om)$ and the cohomology groups $H^{0,2}_J (TM)$ and 
$H^{0,2}_J (M)$ are zero, then $\I_\om$ is a submanifold of $\Jj_\om$ 
in the neighborhood of $J$, with tangent space
\[
T_J \I_\om = \ker \left\{ \db:S\Om^{0,1}_J(TM) \to
  \Om^{0,2}_J(TM) \right\} \subset S\Om^{0,1}_J (TM) = T_J \Jj_\om\,. 
\]
Moreover, the moduli space of infinitesimal compatible deformations of 
$J$ in $\I_\om$ coincides with the moduli
space of infinitesimal deformations of $J$ in $\I$, i.e. it is given by 
$H_J^{0,1}(TM)$ and
\[
T_J \I_\om \cong T_J ((\Diff(M)\cdot J)\cap \I_\om) \oplus H^{0,1}_J (TM)
\,.  
\]
\end{thm}

\begin{remark}
If $H^1(M,\R)=0$, then $cl\Om_J^{0,1}(M) = \db \Om^0 (M,\C) = \Om^0
(M,\C) / \C$ is naturally identified with the complexified Lie algebra
of the symplectomorphism group $\Symp(M,\omega)$.
Moreover, the kernel of $\delta$ can be identified with the vector space of holomorphic vector
fields on $(M,J)$. This nicely agrees with Donaldson's formal picture
where
\[
(\Diff(M)\cdot J) \cap \I_\om = 
(\Symp(M,\omega)^\C \cdot J)\,.
\]
\end{remark}

\subsection{Diffeomorphic compatible complex structures:}

\noindent

\noindent
In this subsection $(M,\omega)$ will be an arbitrary symplectic
manifold, not necessarily of dimension $4$, and $\I_\om$ 
will denote again the space of complex structures on $M$ compatible
with $\omega$.  
We write 
\begin{equation}
\label{diffpressymp}
\Diff_{[\omega]}(M) = \{ \varphi \in \Diff(M) \mid 
\varphi^*([\omega])=[\omega] \in H^{2}(M;\R) \} 
\end{equation}
for the subgroup of diffeomorphisms of $M$ preserving the cohomology
class of the symplectic form. 

Let $J_0\in\I_\om$ denote a fixed complex structure
compatible with $\omega$. Define 
\begin{align}
& U = \{ J \in \I_\om \mid J_0 = \varphi^*J \text{ for
  some } \varphi \in \Diff_{[\omega]}(M) \} \subset \I_\om
\notag \\
& \mathit{\Omega} = \{\eta \in \Omega^2(M) \mid d\eta=0, [\eta]=[\omega] 
\text{ and } \eta \text{ is compatible with } J_0 \} \,. \notag 
\end{align}
$\mathit\Omega$ is a contractible convex subset in $\Omega^2(M)$ and we want to
describe the topology of $U$.

Define a map from $\mathit\Omega$ to the identity component $\Diff_0(M)$ of the
diffeomorphism group of $M$, 
\[ 
\Psi : \mathit\Omega \to \Diff_0(M) \,,
\]
as follows. Given $\eta \in \mathit\Omega$, let $\psi_t \in \Diff_0(M)$ be
the isotopy satisfying $\psi_t^*((1-t)\omega + t\eta) = \omega$ which
is canonically determined by Moser's method and the Riemannian metric
given by $(\om,J_0)$. Then
\[  
\Psi(\eta):=\psi_1\,. 
\]
Note that, if $K = \Iso(\om, J_0)$ denotes the K\"ahler isometry group and $g\in K$ then we have 
$\psi_t(g^*\eta) = g\psi_t(\eta)$ and, in particular,
$\Psi(g^\ast\eta) = g\Psi(\eta)\,,\ \forall\,\eta\in\mathit\Omega$.
Moreover, $\Psi(\eta)^\ast (\eta) = \om \,,\ \forall\,\eta\in\mathit\Omega$.

Denote by $\Hol_{[\om]}(J_0)$ the group of complex automorphisms of
$(M,J_0)$ that preserve the cohomology class $[\om]\in H^2(M;\R)$.
\begin{prop}
The map $\mu: \Symp(M,\omega) \times \mathit\Omega \to U$ defined by 
\[ 
\mu(\phi, \eta) = (\phi^{-1})^*\Psi(\eta)^* J_0 
\]
is a principal $\Hol_{[\om]}(J_0)$-bundle.
\end{prop}
\begin{proof}
Suppose $J \in U$ and let $\varphi \in \Diff_{[\omega]}(M)$ be such
that $\varphi^*(J)=J_0$. Setting $\eta=\varphi^*\omega$, we derive:
\begin{align}
& J_0^* \eta = J_0^*\varphi^*\omega = (d\varphi \circ J_0)^*\omega 
     = (J \circ d\varphi)^*\omega = \varphi^*J^*\omega 
     = \varphi^*\omega =\eta  \notag \\
& \eta(X,J_0X) = \omega((d\varphi) X, (d\varphi) J_0 X) 
     = \omega((d\varphi) X, J (d\varphi) X) > 0 \,. \notag
\end{align}

\noindent
so we conclude that $\eta$ is compatible with $J_0$. 
Since $[\eta]=[\omega]$ we see that $\eta \in \mathit\Omega$.
Moreover, $\varphi\Psi(\eta) \in \Symp(M,\omega)$, since 
$(\varphi \Psi(\eta))^*\omega = \Psi(\eta)^*\eta =\omega$,
and we have
\[ 
\mu(\varphi \Psi(\eta), \eta) = (\Psi(\eta)^{-1} \varphi^{-1})^* 
\Psi(\eta)^*J_0 = (\varphi^{-1})^* J_0 = J \,.
\]
Hence, the map $\mu$ is surjective. 

\medskip
\noindent
Now, given $J\in U$ and an element $(\phi,\eta)\in\mu^{-1}(J)$, consider $\varphi = \phi \Psi(\eta)^{-1}$. Then $\varphi \in \Diff_{[\om]} (M)$ and
\[
\mu (\phi,\eta) = (\phi^{-1})^* \Psi(\eta)^* J_0 = J \Rightarrow
J_0 = \varphi^* J\,.
\]
Conversely, given $\varphi\in \Diff_{[\om]} (M)$ such that $J_0 =
\varphi^* J$, we have that $(\varphi \Psi(\eta), \eta) \in \mu^{-1}
(J)$, where $\eta = \varphi^*\om$. Hence
\[ 
\mu^{-1}(J) \cong \{ \varphi \in  \Diff_{[\om]} (M)\ \ | \ \ 
J_0 = \varphi^*J \} \,,
\]
i.e. the fibers of the map $\mu$ are torsors on the group
$\Hol_{[\om]} (J_0)$. In fact, defining a right action of $\Hol_{[\om]} (J_0)$ on 
$\Symp(M,\omega) \times \mathit\Omega$ by
\[ 
(\phi,\eta) \cdot \varphi = 
(\phi\Psi(\eta)^{-1}\varphi\Psi(\varphi^*(\eta)),\varphi^*(\eta))\,, 
\]
we see that the fibers of the map $\mu$ are free orbits of this
action.
\end{proof}

\begin{cor} \label{cor:intorbit}
If $J_0 \in \I_\om$ is such that the inclusion
$\Iso(\om,J_0) \hookrightarrow \Hol_{[\om]} (J_0)$ is a weak homotopy 
equivalence, then the inclusion of the $\Symp(M,\om)$-orbit of $J_0$ in
$U$, i.e.
\[ 
\Symp(M,\omega)/\Iso(\om,J_0) \hookrightarrow U \,,
\]
is also a weak homotopy equivalence.
\end{cor}
\begin{proof}
Indeed, we have 
\[
\Symp(M,\om) \stackrel{\sim}{\hookrightarrow} \Symp(M,\om) \times \{\om\}
\subset \Symp(M,\om) \times \mathit\Omega \stackrel{\mu}{\rightarrow} U\,,
\]
which induces
\[ 
\Symp(M,\om)/\Iso(\om,J_0) \stackrel{\sim}{\hookrightarrow} 
\left(\Symp(M,\om) \times \mathit\Omega\right)/\Iso(\om,J_0)  
\stackrel{\overline{\mu}}{\rightarrow} U\,,
\]
where the fiber of $\overline{\mu}$ is weakly contractible by assumption.
\end{proof}

\begin{remark} \label{rmk:intorbit}
According to Calabi~\cite{C1}, a source of examples for the previous 
corollary are $J_0$'s determining extremal K\"ahler metrics, at least 
on manifolds where $\Hol(J_0)$ and $\Iso(\om,J_0)$ are both connected.
\end{remark}

\subsection{Transversality:}

\noindent

\noindent
In this subsection we study transversality properties of certain
strata of compatible almost complex structures on a compact symplectic
manifold $(M,\om)$, with respect to the inclusion $\I_\om
\subset \Jj_\om$. These strata are characterized by the existence of
certain pseudo-holomorphic curves.

Recall that if $(M,J)$ is an almost complex manifold and $(\Sigma, j)$
is a Riemann surface, a map $u\in\Map (\Sigma, M)$ is called a
$J$-holomorphic curve if
\[
J \circ du = du \circ j\,,
\]
and a simple $J$-holomorphic curve if, in addition, it is not multiply
covered (see~\cite{MS}).

Given a compact symplectic manifold $(M,\om)$, a compact Riemann
surface $(\Sigma, j)$ and a homology class $A\in H_2 (M,\Z)$, consider
the space
\begin{align}
\Mm (A, \Sigma) = \{ (u,J)\in\Map(\Si,M)\times\Jj_\om\,:\  
& \text{$u$ is a simple $J$-holomorphic curve} \notag \\
& \text{with $u_\ast([\Si])=A$}\}\,. \notag
\end{align}
Denote by $\Om^{0,1}(\Si, TM)$ the vector bundle over $\Map(\Si,M)
\times \Jj_\om$ whose fiber over $(u,J)$ is
\[
\Om^{0,1}_J (\Si, u^\ast TM) \equiv \ \text{$J$ anti-linear $1$-forms
  on $\Si$ with values in $u^\ast TM$.}
\]
This vector bundle has a natural section, denoted by $\db$, given at
$(u,J)\in \Map(\Si,M)\times\Jj_\om$ by
\[
\db_{(u,J)} \equiv \db_J(u) \equiv \frac{1}{2} (du + J\circ du\circ j)
\in \Om^{0,1}_J (\Si, u^\ast TM)\,.
\]
The space $\Mm (A, \Sigma)$ is the zero set of this section.

The tangent space to $\Map(\Si,M)\times\Jj_\om$ at $(u,J)$ is given by
\[
T_{(u,J)}(\Map(\Si,M)\times\Jj_\om) = T_u \Map(\Si, M) \oplus T_J
\Jj_\om = \Om^0 (\Si, u^\ast TM) \oplus S\Om^{0,1}_J (TM)\,.
\]
At $(u,J)\in\db^{-1}(0) = \Mm(A, \Si)$, the vertical component of the
derivative of the section $\db$,
\[
D(\db)_{(u,J)} : \Om^0 (\Si, u^\ast TM) \oplus S\Om^{0,1}_J (TM) \to
\Om^{0,1}_J (\Si, u^\ast TM)\,,
\]
is surjective and given by
\[
D(\db)_{(u,J)}(\xi, \alpha) = \db_J (\xi) + u^\ast (\alpha)\,,
\]
where $\db_J(\xi)$ is given as in Definition~\ref{defn:db0TM} of
Appendix~\ref{app:A} (see~\cite{MS}, Remark 3.1.2 and Proposition
3.2.1). Hence, the following holds.
\begin{prop} \label{prop:mspace}
$\Mm (A,\Si)$ is an infinite-dimensional submanifold of $\Map (\Si, M)
\times \Jj_\om$, with tangent space given by
\[
T_{(u,J)} \Mm(A, \Si) = \{ (\xi,\alpha) \in 
\Om^0 (\Si, u^\ast TM) \oplus S\Om^{0,1}_J (TM)\,:\ 
\db_J (\xi) + u^\ast (\alpha) = 0 \}\,.
\]
\end{prop}
\noindent
The image of the projection
\[
\pi : \Mm (A,\Si)  \to \Jj_\om \;, \quad  (u,J)  \mapsto J
 \]
defines a subset $U_A\subset \Jj_\om$ characterized by
\begin{align}
J\in U_A \Leftrightarrow \ & \text{$A\in H_2(M,\Z)$ can be represented
  by a simple} \notag \\
& \text{$J$-holomorphic curve with domain $(\Si, j)$.} \notag
\end{align}
We now want to find conditions ensuring that the image $U_A$ of $\pi$ 
is transversal to
$\I_\om \subset \Jj_\om$ at points $(u,J) \in \Mm (A, \Si)$
with $J\in\I_\om$. It follows from the previous proposition
that
\[
\pi_\ast (T_{(u,J)}\Mm(A,\Si)) = \{ \alpha \in S\Om_J^{0,1}(TM)\,:\ 
[u^\ast \alpha] = 0 \in H_J^{0,1} (\Si, u^\ast TM)\}\,.
\]
Assume that the restriction map
\[
u^\ast: clS\Om_J^{0,1} (TM) \to H_J^{0,1} (\Si, u^\ast TM)
\]
is surjective. Then, given any $\gamma \in T_J \Jj_\om = S\Om_J^{0,1}
(TM)$, there exists $\beta\in T_J \I_\om =
clS\Om_J^{0,1}(TM)$ such that $(\gamma - \beta) \in \pi_\ast
(T_{(u,J)}\Mm(A,\Si))$. Combining this with the reasoning leading to
Theorem~\ref{thm:intsubacs}, one gets the following result.
\begin{thm} \label{thm:transversal}
Let $(M, \om, J\in\I_\om)$ be a K\"ahler $4$-manifold such
that the cohomology groups $H_J^{0,2}(M)$ and $H_J^{0,2} (TM)$ are
zero. Suppose that $(u,J)\in\Mm(A, \Si)$ is such that 
$u^\ast : H_J^{0,1}(TM) \to H_J^{0,1}(u^\ast(TM))$ is an isomorphism. 
Then $\pi:\Mm (A, \Si) \to \Jj_\om$ is transversal at $(u,J)$ to 
$\I_\om \subset \Jj_\om$ and the infinitesimal complement to the image 
$U_A$ of $\pi$ in a neighborhood of $J$ can be identified with the 
moduli space of infinitesimal deformations $H_J^{0,1}(TM)$.
\end{thm}

\section{Compatible Complex Structures on Rational Ruled Surfaces}

\noindent
In this section we describe the stratification on the space
$\Jj_{\omega}$ of compatible almost complex structures on a rational
ruled surface, previously studied in~\cite{Ab,AM, McD1}, and use
results from the previous section to show that it induces an analogous
stratification on the corresponding space $\I_{\omega}$ of
compatible complex structures. We then conclude with the proof of
theorem~\ref{Main1}.

In the rest of the paper we will denote by $(M,\omega_\lambda)$ 
the symplectic manifolds
$S^2 \times S^2$, with the split symplectic form
\begin{equation}
\label{untwist}
\omega_\lambda = \la\sigma \oplus \sigma,
\end{equation} 
with $1\leq\la\in\R$, as well as the nontrivial bundle $S^2\tilde{\times} S^2$ with a symplectic form 
\begin{equation}
\label{twisted}
\omega_\la \quad  \text{ satisfying } \quad  [\om_\la](E) = \la, \quad [\om_\la](F) = 1,
\end{equation}
with $0<\la \in \R$, where $F$ is the homology class of a fiber and $E$ is the homology class of the exceptional divisor under the natural identification of $S^2 \tilde{\times} S^2$ with $\CP^2 \sharp \overline{\CP^2}$.

We will refer to the case of $M=S^2\times S^2$ as the \emph{untwisted case} and to $M=S^2\tilde{\times}
S^2$ as the \emph{twisted case}. We identify $H_2(S^2\times S^2;\Z)$ with $\Z \times \Z$ in the standard way and  $H_2(S^2\tilde{\times} S^2; \Z)$ with $\Z\times \Z$ via
\begin{equation}
\label{basistwisted}
H_2(S^2\tilde{\times} S^2; \Z) \ni mE+nF \mapsto (m,n)\in
\Z\times\Z\,. 
\end{equation}
Finally, we will write $\Jj_\la$ for the contractible space of almost complex structures on $M$ compatible with $\omega_\la$ and $G_\la$ for the symplectomorphism group of $(M,\omega_\lambda)$.

\subsection{Compatible almost complex structures:}

\noindent
The following theorem plays a fundamental role in the results obtained
in~\cite{Ab} and~\cite{AM} regarding the topology of $G_\la$. 
It will also play a fundamental role here. The most technical point, 
listed as (v) in the statement, was proved in~\cite[Theorem 1.2]{McD1} using 
gluing techniques for pseudo-holomorphic spheres.

\begin{thm} \label{thm:acs}
Let $(M,\om_\la)$ be one of the symplectic manifolds defined above (see \eqref{untwist} and \eqref{twisted}).
There is a stratification of the contractible space $\Jj_\la$ of compatible almost
complex structures of the form
\[ \Jj_\la = U_0 \sqcup U_1 \sqcup \cdots \sqcup U_\ell\,,\]
with $\ell\in\N_0$ such that $\ell < \la \leq \ell +1$ satisfying:
\begin{enumerate}
\item[(i)]
\begin{align}
U_k \equiv \{J\in \Jj_\la &:\,\text{ $(1,-k)\in H_2 (M;\Z)$ is
    represented} \notag \\
& \quad\text{by a $J$-holomorphic sphere}\}\,. \notag
\end{align}
\item[(ii)] $U_0$ is open and dense in $\Jj_\la$. For $k\geq 1$, $U_k$
  has codimension $4k-2$ in $\Jj_\la$ in the untwisted case and codimension $4k$ in the twisted case.
\item[(iii)] $\overline{U_k} = U_k \sqcup U_{k+1} \sqcup \cdots \sqcup U_\ell$.
\item[(iv)] Given a compatible almost complex structure $J_k \in U_k$ with isometry group
\[
\Iso(\om_\la,J_k) \cong 
\begin{cases}
\Z/2 \ltimes (SO(3)\times SO(3))\,, &\text{if $\la = 1$ and $k=0$, in the untwisted case,} \\
SO(3)\times SO(3)\,, &\text{if $\la > 1$ and $k=0$, in the untwisted case,} \\
S^1\times SO(3)\,, &\text{if $k\geq 1$, in the untwisted case,} \\
U(2)\,, & \text{ in the twisted case,}
\end{cases} 
\]
then the inclusion
\[ 
\left(G_\la / Iso(\om_\la,J_k) \right)  \lra U_k \;, \quad [\psi]  \longmapsto 
\psi_\ast (J_k) 
\]
is a weak homotopy equivalence.
\item[(v)] Each $U_k$ has a tubular neighborhood $NU_k \subset
  \Jj_\la$ which fibers over $U_k$ as a ball bundle.
\end{enumerate}
\end{thm}

\begin{remark}
In~\cite{AM}, Section 2, such $(S^2\times S^2, \om_\la, J_k)$ and
$(S^2\tilde{\times} S^2, \om_\la, J_k)$ were explicitly constructed
as K\"ahler reductions of $\C^4$. This is reviewed below in 
Section \ref{deformations}.
\end{remark}

\subsection{Compatible complex structures:}

\noindent

\noindent
Our goal now is to show that an analogous theorem holds for the space
$\I_\la\subset\Jj_\la$ of compatible integrable complex
structures on $(M,\om_\la)$. 

For each $k\in\left\{0,1,\ldots,\ell\right\}$, with $\ell\in\N_0$ 
such that $\ell < \la \leq \ell +1$, define
\begin{align}
V_k \equiv U_k \cap \I_\la =
\{J\in\I_\la&:\,\text{$(1,-k)\in H_2 (M;\Z)$ is
    represented} \notag \\
& \quad\text{ by a $J$-holomorphic sphere}\}\,. \notag
\end{align}
It follows from standard complex geometry (see~\cite[Proposition V.4.3(i)]{BPV}
or~\cite{Ca}) that any $J \in V_k$ is complex isomorphic to
the Hirzebruch surface
\[
F_{n} = \CP (\Oo\oplus\Oo(-n))\,,\ \text{with $n=2k$ if $M =
  S^2\times S^2$, and $n=2k+1$ if $M = S^2\tilde{\times} S^2$,}
\]
by a diffeomorphism of $M$ that acts as the identity in homology
(here, $\Oo(-1)$ denotes the tautological line bundle over $\CP^1$
and $\CP(E)$ the projectivization of a vector bundle $E$). 
Moreover, Calabi proved in~\cite{C1} that there is a complex structure
$J_k \in V_k$, unique up to the action of $G_\la$, for which 
$g_{\la,k} \equiv \om_\la (\cdot, J_k \cdot )$ is an extremal K\"ahler
metric, with K\"ahler isometry group $K_k \equiv \Iso (\om_\la, J_k)$
as in Theorem \ref{thm:acs}. These two facts
together imply that
\[
V_k = \{ J \in \I_\la \mid J_k = \varphi^*J \text{ for
  some } \varphi \in \Diff_{[\omega_\la]}(M) \} \subset \I_\la.
\]
\begin{thm} \label{thm:intorbit}
The inclusions
\[
G_\la / K_k \hookrightarrow V_k \hookrightarrow U_k
\]
are weak homotopy equivalences.
\end{thm}
\begin{proof}
Corollary~\ref{cor:intorbit} says that the left map is a 
weak homotopy equivalence while Theorem~\ref{thm:acs}(iv)
says that the composite is a weak
homotopy equivalence. It follows that the map $V_k \hookrightarrow U_k$ is also
a weak homotopy equivalence.
\end{proof}

Any Hirzebruch surface $F_{n}$ satisfies $H^{0,2}(F_{n}) =
H^{0,2}(TF_{n}) = 0$ \cite{Ko}. Hence, it follows from
Theorem~\ref{thm:intsubacs} that $\I_\la$ is an infinite
dimensional submanifold of $\Jj_\la$. Since $U_0$ is open in
$\Jj_\la$, we have that $V_0$ is also open in
$\I_\la$.

For each $k\in\{1, \ldots, \ell\}$, consider the space
\begin{align}
\Mm_k = \{ (u,J)\in\Map(S^2, M)\times\Jj_\la\,:\  
& \text{$u$ is a simple $J$-holomorphic sphere with} \notag \\
& \text{$u_\ast([S^2])= (1,-k)\in H_2 (M; \Z)$}\}\,, \notag 
\end{align}
which by Proposition~\ref{prop:mspace} is an infinite dimensional
submanifold of $\Map(S^2, M)\times\Jj_\la$. Positivity of
intersections and the adjunction inequality for pseudo-holomorphic
curves in almost complex $4$-manifolds (see Theorems 2.6.3 and 2.6.4
in~\cite{MS}) imply that: 
\begin{itemize}
\item[-] if $(u,J)\in\Mm_k$ then the simple $J$-holomorphic map $u:S^2
  \to M$ is an embedding;
\item[-] if $(u_1,J), (u_2, J) \in \Mm_k$ then the simple
  $J$-holomorphic maps $u_1, u_2:S^2 \to M$ have exactly
  the same image, i.e. they differ only by an holomorphic
  reparametrization of $S^2$ given by an element in $PSL(2,\C)$.
\end{itemize}
This means that the projection
\[ 
\pi : \Mm_k  \to U_k \subset \Jj_\la \;, \quad (u,J)  \mapsto J 
\]
is a principal $PSL(2,\C)$-bundle map over $U_k$. 

\begin{prop}
For any $(u,J)\in\Mm_k$, the map $u^\ast: H^{0,1}(TF_{n}) \to 
H^{0,1}(u^\ast (TF_{n}))$ is an isomorphism, where $n=2k$ if 
$M = S^2\times S^2$ and $n=2k+1$ if $M = S^2\tilde{\times} S^2$. 
\end{prop}
\begin{proof}
Recall that $F_{n} = \CP(\Oo \oplus \Oo(-n))$. The inclusion map 
$u : \CP^1 \rightarrow F_{n}$ corresponds to the zero section 
$\CP(0\oplus \Oo(-n))$. Let $v : \CP^1 \rightarrow F_{n}$ denote the 
section at infinity: $\CP(\infty \oplus \Oo(-n))$. Let 
$i: Z(0) \subset F_{n}$ be the complement of $u$, and 
$j: Z(\infty) \subset F_{n}$ be the complement of $v$. 
Notice that the space $Z(0)$ is equivalent to the total space of the 
bundle $\Oo(n)$ over $v$, and $Z(\infty)$ is equivalent to $\Oo(-n)$ 
over $u$. We have a short exact sequence of sheaves of 
$\Oo_{F_{n}}$-modules:
\[ 
0 \rightarrow \Oo_{F_{n}} \rightarrow j_*j^* \Oo_{F_{n}} 
\rightarrow v_* v^! \Oo_{F_{n}} \rightarrow 0 
\]
where $v^!$ is a functor from $\Oo_{F_{n}}$-modules to 
$\Oo_{\CP^1}$-modules such that the stalk of $v^!\mathcal{S}$ at 
$x \in \CP^1$ (identified with the image of $v$) is given by the quotient:
\[ 
\mathcal{S}_x \rightarrow ( j_*j^* \mathcal{S})_x \rightarrow (v^! 
\mathcal{S})_x \rightarrow 0 
\]
Hence, $v^!\mathcal{S}$ may be seen as the higher residues along the normal. 
Identifying $Z(0)$ with $\Oo(n)$, it follows that $v^!\Oo_{F_{n}} = \mbox{Sym}_+\Oo(n)$, 
where $\mbox{Sym}_+\mathcal{S}$ stands for the augmentation ideal in 
the symmetric algebra on $\mathcal{S}$. We may now tensor the above 
short exact sequence with $TF_{n}$ to get:
\[ 
0 \rightarrow \Oo(TF_{n}) \rightarrow j_*j^* \Oo(TF_{n}) 
\rightarrow v_* v^! \Oo(TF_{n}) \rightarrow 0 
\]
In cohomology, we get the following exact sequence:
\[ 
\ldots \rightarrow H^1(TF_{n}) \rightarrow H^1(Z(\infty),
j^*TF_{n}) \rightarrow H^1(\CP^1,v^!TF_{n}) \rightarrow \ldots
\]
Since $v^!TF_{n}$ is the bundle 
$(\Oo(2) \oplus \Oo(n)) \otimes \mbox{Sym}_+\Oo(n)$, 
the last term in the above sequence is trivial. Therefore the restriction map:
\[ j^*: H^1(TF_{n}) \rightarrow H^1(Z(\infty), j^*TF_{n})\]
is an epimorphism. Now recall that $Z(\infty)$ is the total space of the bundle $\Oo(-n)$ over $u$. Let $\pi : Z(\infty) \rightarrow \CP^1$ be the projection map for this bundle. We get an isomorphism: 
\[ \pi_* : H^*(Z(\infty),j^*TF_n) \rightarrow H^*(\CP^1,\pi_*j^*TF_n) = H^*(\CP^1,u^*TF_n\otimes \mbox{Sym} \, \Oo(n)) \]
The long exact sequence in cohomology now gives us an epimorphism:
\[ u^* : H^1(Z(\infty), j^*TF_{n}) \rightarrow H^1(\CP^1, u^*TF_n\otimes \mbox{Sym} \, \Oo(n)) \rightarrow  H^1(\CP^1, u^*TF_{n}) \]
Composing the above sequence of epimorphisms, we see that the required map $u^*$ in the statement of the proposition is an epimorphism. Finally, notice that $u^*TF_{n} = \Oo(2) \oplus \Oo(-n)$ , hence the dimension of  $H^1(\CP^1,u^*TF_{n})$ and $H^1(TF_{n})$ both equal $n-1$ 
(see for example \cite[Example 6.2(b)(4), p.309]{Ko}. The proof follows.
\end{proof}

\begin{remark}
See also \cite[Lemma A.8]{ALP} for a generalization of the previous Proposition to blow-ups of a rational ruled surface.
\end{remark}

\noindent
It follows from Theorem~\ref{thm:transversal} that, for $k\in\{1, \ldots,
\ell\}$, each strata $U_k \subset \Jj_\la$ is transversal to
$\I_\la \subset \Jj_\la$. Hence, the stratification of
$\Jj_\la$ induces by intersection a stratification of $\I_\la$ of the form
\[
\I_\la = V_0 \sqcup V_1 \sqcup \cdots 
\sqcup V_\ell\,,
\]
which satisfies the direct analogues of items (i), (ii), (iii) and
(iv) in Theorem~\ref{thm:acs}.

Since each stratum $V_k$, $k\in\{1, \ldots, \ell\}$, is the
transversal intersection of $U_k$ and $\I_\la$, the tubular
neighborhood $NU_k\subset \Jj_\la$ of Theorem~\ref{thm:acs}-(v) gives
rise to a tubular neighborhood $NV_k \equiv NU_k \cap
\I_\la$ of $V_k$ in $\I_\la$, which
fibers over $V_k$ as a ball bundle. By
Theorem~\ref{thm:transversal}, each of these balls can be identified
with a neighborhood of zero in
\[
H^{0,1}(TF_{n}) \cong \C^{n-1}\quad\text{\cite[Example 6.2(b)(4), p.309]{Ko}.}
\]

\subsection{Proof of Theorem~\ref{Main1}:}

\noindent

\noindent
It follows from the results stated in the previous two subsections
that the inclusion $\I_\la \hookrightarrow \Jj_\la$ is
transversal to the stratification
\[
\Jj_\la = U_0 \sqcup U_1 \sqcup \cdots \sqcup U_\ell\,.
\]
By Theorem \ref{thm:intorbit}, the induced stratification
\[
\I_\la = V_0 \sqcup V_1 \sqcup \cdots 
\sqcup V_\ell
\]
is such that the inclusions
\[
V_k \hookrightarrow U_k \quad\text{and}\quad
NV_k \setminus V_k \hookrightarrow NU_k \setminus U_k
\]
are weak homotopy equivalences. Writing $U_{0i} = U_0 \sqcup \cdots \sqcup U_i$ and 
$V_{0i} = V_0 \sqcup \cdots \sqcup V_i$, assume inductively that the inclusion
\[ U_{0(i-1)} \to V_{0(i-1)} \]
is a weak equivalence. Then we have a map of excisive triads
\[ (V_{0i}, NV_i, V_{0(i-1)}) \to (U_{0i}, NU_i, U_{0(i-1)}) \] 
which is a weak equivalence when restricted to $NV_i, V_{0(i-1)}$ and their 
intersection $NV_i\setminus V_i$. It follows (see for instance \cite[p. 80]{May}) that 
\[V_{0i} \to U_{0i}\]
is a weak equivalence and the result now follows by induction.

Since $\Jj_\la$ is contractible, saying that $\I_\la \subset \Jj_\la$ is a weak equivalence 
is of course equivalent to the statement that $\I_\la$ is weakly contractible.

\begin{remark}
\label{remarktame}
Note that the arguments above also apply to the space 
of complex structures tamed by a symplectic form. Thus Theorem~\ref{Main1}
still holds with the word compatible replaced by tame.
\end{remark}
 
\section{Deformation of Hirzebruch Surfaces}

\label{deformations}

\noindent
We will use the notation fixed in the introduction of the previous section.
In that section, we described a stratification of the space $\I_\la$ of compatible integrable complex structures on $(M,\omega_\lambda)$ which is equivariant with respect to the action of $G_\la$.
 
Our aim in this section is to compute the representations of the isotropy
groups $\Iso(\omega_\lambda,J_k)$ of Theorem \ref{thm:acs} on the links of the strata $V_k$ containing $J_k$. According to Theorem \ref{thm:transversal}, the links can be identified with $H^{0,1}(TF_n)$ where $n=2k$ in the untwisted case and $n=2k+1$ in the twisted case (cf. discussion in the previous section). Atiyah and Bott's fixed point theorem \cite[Theorem II.4.12]{AB2} reduces this calculation to a calculation of the isotropy representations of the maximal torus of $\Iso(\omega_\lambda,J_k)$ on the tangent spaces to its (four) fixed points in $M$. In order to obtain these, we will first describe a construction of $(\omega_\lambda,J_k)$ by K\"ahler reduction (which appears for example in \cite[Section 2.3]{AM}). This construction will also give us a hold on the K\"ahler isometry groups in Theorem \ref{thm:acs} and allow us to relate the standard basis for their maximal tori to the one arising naturally from the K\"ahler reduction procedure. All of this will be necessary
  in the next section when we compute the
cohomology of the classifying space of the symplectomorphism groups.

\subsection{K\"ahler reduction:}

\noindent

\noindent
K\"ahler reduction of $\C^4$ by the action of the $2$-torus $T^2_n$ acting via
\[ (s,t) \cdot (z_1,\ldots,z_4) = (s^ntz_1,tz_2,sz_3,sz_4) \]
at the values
\[
\begin{cases}
(\la + \tfrac{n}{2},1) & \text{ if } n \text{ is even,}\\
(\la + \tfrac{n+1}{2},1) & \text{ if } n \text{ is odd,}
\end{cases}
\]
produces K\"ahler manifolds which are symplectomorphic to $(M,\omega_\lambda)$ where $M=S^2\times S^2$
if $n$ is even and $M=S^2\tilde \times S^2$ if $n$ is odd (see \cite[Section 2.3]{AM} for a detailed explanation of this). As complex manifolds these reductions are isomorphic to Hirzebruch surfaces $F_n$ with projection onto the base $\CP^1$ given by $[(z_1,z_2,z_3,z_4)] \mapsto [z_3\colon z_4]$.

The inclusion of the torus $T^2_n$ in the standard torus $T^4 \subset U(4)$ is given by the matrix
\[ \begin{bmatrix} n & 1 \\ 0 & 1 \\ 1 & 0 \\ 1 & 0 \end{bmatrix} \]
The connected component of the K\"ahler isometry group of $F_n$ is 
\[ K(n)=Z_{U(4)}(T^2_n)/T^2_n = \begin{cases} 
(U(2) \times U(2))/T^2_0 & \text{ if  } n=0, \\
(T^2 \times U(2))/T^2_n & \text{ if }  n\geq 1, 
\end{cases}
\]
and so
\begin{equation}
\label{isomgroups}
K(n) = \begin{cases} 
SO(3) \times SO(3) & \text{ if } n=0 \text{ in the untwisted case,} \\
S^1 \times SO(3) & \text{ if } n>0 \text{ is even,} \\
U(2) & \text{ otherwise.}
\end{cases}
\end{equation}
The maximal torus $\T$ of $K(n)$ is $T^4/T^2_n$. We will identify it with $T^2$ via the composite
\begin{equation}
\label{momentbasis}
T^2 \xrightarrow{\phi} T^4 \to T^4/T^2_n 
\end{equation}
with
\[ \phi = \left[\begin{matrix} 0 & 0 \\ 1 & 0 \\ 0 & 1 \\ 0 & 0 \end{matrix}\right] \]
We will refer to this basis as \emph{the moment map basis} for the maximal torus $\T$.
We will need to relate this basis for $\T$ with the standard basis arising from the identification \eqref{isomgroups}. First note that
\begin{equation}
\label{torusso3}
SO(3)=PU(2)=U(2)/\Delta(S^1)
\end{equation}
and so there is a natural quotient map $U(2) \to SO(3)$. We take the standard torus inside $SO(3)$ to be the image of $U(1)\times 1 \subset U(2)$ under this quotient, and we make the obvious choice for basis of $T^2\subset U(2)$.

For $n>0$ we have $K(n) = U(2)_{z_3,z_4}/(\Z/n)$ which can be naturally identified with the groups
described in \eqref{isomgroups} using the double cover $S^1 \times SU(2) \to U(2)$ and the canonical isomorphism $S^1/(\Z/k) \simeq S^1$. It is not hard to check that with respect to the standard bases 
described above and these identifications, the $n$-fold covering map $U(2)_{z_3,z_4} \to K(n)$ is given on maximal tori by the matrix 
\begin{equation}
\label{coverings}
\begin{bmatrix} \frac{n+1}2 & \frac{n-1}2 \\ \frac{n-1}2 & \frac{n+1}2 \end{bmatrix} \text{ for } n \text{ odd, and } \quad  \begin{bmatrix} \frac{n}{2} & \frac{n}{2} \\  1 & -1 \end{bmatrix} 
\text{ for } n>0 \text{ even.}
\end{equation}

\begin{lemma}
\label{changeofbasis}
The changes of basis from the moment map bases \eqref{momentbasis} to the standard bases for the 
maximal tori of $K(n)$ are given by the matrices
\[ \begin{bmatrix}  1 & \frac{n+1}2 \\ 1 & \frac{n-1}2 \end{bmatrix} \text{ for } n \text{ odd, } \quad \begin{bmatrix} 1 & \frac n 2 \\ 0 & 1 \end{bmatrix} \text{ for } n>0 \text{ even, and}
\quad \begin{bmatrix} -1 & 0 \\ 0 & 1 \end{bmatrix} \text{ for } n=0. \]
\end{lemma}
\begin{proof}
Suppose $n>0$ is odd. The image of the circle $(0,n,0,0)$ in $K(n)$ equals the image of $(0,0,1,1)$. This is
the circle $(1,1)$ in $U(2)_{z_3,z_4}$ which according to \eqref{coverings} is taken to the 
circle $(n,n)$ in $K(n)$. Thus the first element of the moment map basis \eqref{momentbasis} is
mapped to the standard diagonal in $K(n)$ and we get the first column of the first matrix in the 
statement. The other computations are similar.
\end{proof}

\subsection{Isotropy representations:}

\noindent

\noindent
We now want to compute the representation of $K(n)$ on $H^{0,1}(TF_n)$ which,
according to Section 3, is the isotropy representation of $\Iso(\omega_\lambda),J_k)$
on the link of the stratum $V_k \subset \I_\la$.

For this we will use the $K(n)$-equivariant elliptic complex
\begin{equation}
\label{complex}
\Omega^0(TF_n) \xrightarrow{\db} \Omega^{0,1}(TF_n) \xrightarrow{\db} 
\Omega^{0,2}(TF_n). 
\end{equation}
Since $H^{0,2}(TF_n) = 0$ (see \cite[Example 6.2(b)(4), p.309]{Ko}), the 
index of this complex is the virtual representation
\[ H^0(TF_n) - H^{0,1}(TF_n) \]
of $K(n)$. The Atiyah-Bott fixed point theorem \cite[Theorem 4.12]{AB2} gives a 
formula for the character of this representation in terms of the characters of the isotropy representation
of the maximal torus $\T$ of $K(n)$ on the fibers of $TF_n$ over the $\T$-fixed points.
We start by determining the latter using the construction of $F_n$ by K\"ahler reduction given
in the previous subsection. 

\medskip
\noindent
Using the basis \eqref{momentbasis}, the moment map for the action of $\T$ on $F_n$ is given by the expression
\[ [(z_1,\ldots,z_4)] \mapsto \frac 1 2 (|z_2|^2,|z_3|^2) \in \R^2=\Lie(\T)^*\]
and so its image is the moment polygon with vertices $A=(0,0), B=(1,0), C=(1,\mu)$ and $D=(0,\mu-n)$ in 
$\R^2$ where
\[ \mu = 
\begin{cases} 
\la + \frac{n}{2} & \text{ if } n \text{ is even,} \\
\la + \frac{n+1}{2} & \text{ if } n \text{ is odd.}
\end{cases}
\]
The vertices of the moment polygon are the images of the four $\T$-fixed points
which we also denote by $A,B,C,D$. The weights of the isotropy representation of the torus on the tangent space at a fixed point are determined by the primitive vectors in $\Hom(\T,S^1)=\Z^2 \subset \R^2=\Lie(\T)^*$ along the edges of the moment polygon which meet at the corresponding vertex. If we write $x$ for the weight determined by the vector $(1,0)$ and $y$ for the weight determined by the vector $(0,1)$ and write the weights multiplicatively, we see that the weights $w_i, i=1,2$ of the two dimensional isotropy representation at the fixed points are given by the following table:
\begin{equation}
\label{weighttable}
\begin{array}{|c|c|c|c|c|} \hline
  \quad  & \quad   A \quad  & \quad B \quad & \quad C \quad & \quad D
\\ 
\hline
w_1   &     x      &    1/x    &   1/(xy^n)        &  xy^n
\\ 
\hline 
w_2   &     y      &   y      &   1/y    &  1/y
\\ 
\hline 
\end{array}
\end{equation}
The character of the isotropy representation of $\T$ on the tangent space at a fixed point is the sum of the two weights in the corresponding column. We can now prove the main result of this section. 
\begin{thm} \label{representation}
Let $n>1$. The representation of $K(n)$ on $H^{0,1}(TF_n)$ is given by 
\[ 
\begin{cases}
\Det^{-\frac{n-3}{2}}\otimes \Sym^{n-2}(\C^2) & \text{ if } n \text{ is odd,}\\
\Det \otimes \Sym^{n-2}(\C^2) & \text{ if } n \text{ is even.} 
\end{cases}
\] 
where $\Det$ denotes the determinant representation of $U(2)$ if $n$ is odd and the standard representation of the $S^1$ factor if $n$ is even; $\Sym^k(\C^2)$ denotes the $k$-th symmetric power of the defining representation of $U(2)$ if $n$ is odd and the irreducible $(k+1)$-dimensional representation of $SO(3)$ if $n$ is even.
\end{thm}
\begin{proof}
By \cite[Theorem II.4.12]{AB2}, the topological index $I(n)$ of \eqref{complex}
is given by the following character\footnote{Note that $TX$ denotes the \emph{cotangent} bundle of $X$ in \cite{AB2}.} of $\T$:
\[ 
I(n) = \sum \frac{w_1+w_2}{\left(1-\frac{1}{w_1}\right)\left(1-\frac{1}{w_2}\right)} = \sum \frac{w_1w_2(w_1+w_2)}{(1-w_1)(1-w_2)},
\]
where the sum is taken over the four fixed points of $\T$ and the weights are given by
\eqref{weighttable}. Writing this in terms of the weights $x$,$y$ yields
\[ 
I(n)= \frac{xy(x+y)}{(1-x)(1-y)} +
\frac{y/x(1/x+y)}{(1-1/x)(1-y)} +
\frac{1/(xy^{n+1})(1/(xy^n)+1/y)}{(1-1/(xy^n))(1-1/y)} + 
\frac{xy^{n-1}(xy^n+1/y)}{(1-xy^n)(1-1/y)} 
\]
and writing this in linearly independent monomials we obtain
\begin{align}
& I(0) = 2 + y + \frac{1}{y} + x + \frac{1}{x}, \quad \quad
\quad 
I(1) = 2 + y + \frac{1}{y} + \frac{1}{xy} (1 + y),    \notag \\
& I(n) = 2 + y + \frac{1}{y} + \frac{1}{xy^n}(1+y+\ldots+y^n) - 
xy(1+y+\ldots+y^{n-2}), \quad \text{ for }n>1. \notag 
\end{align}
The negative terms are the negative of the character of the representation
$H^{0,1}(TF_n)$, which thus has complex dimension $n-1$ for $n>0$, in accordance with the
fact that $\dim_\C H^{0,1}(TF_n) = n-1$ \cite[Example 6.2(b)(4), p.309]{Ko}. 
The positive terms are the character of a representation of dimension $n+5$, namely the restriction
to $K(n)$ of the adjoint representation of the complex automorphism group of $F_n$ on its Lie
algebra $H^0(TF_n)$.

\medskip
\noindent
Using Lemma \ref{changeofbasis} we can write the character of the representation $H^{0,1}(TF_n)$ in terms of the standard weights $a$ and $b$ of $\T \subset K(n)$, dual to the standard basis of the maximal torus. We have 
\[ a = \begin{cases}   xy^{\frac{n+1}{2}} & \text{ if } n>0 \text{ is odd,}\\
xy^{n/2} & \text{ if } n>0 \text{ is even,} \end{cases} \quad \quad
b = \begin{cases}   y^{\frac{n-1}2}x & \text{ if } n>0 \text{ is odd,}\\
y & \text{ if } n>0 \text{ is even.} \end{cases}\]
Hence for $n>0$, the character of the representation of $K(n)$ on $H^{0,1}(TF_n)$ is given by
\begin{equation}
\label{weightsstand}
\begin{cases}
(ab)^{-\frac{n-3}{2}}(b^{n-2} + ab^{n-1} + \ldots + a^{n-2}) & \text{ if } n>1 \text{ is odd,} \\
a(b^{1-\frac n 2} + \ldots + 1 + \ldots + b^{\frac n 2 -1}) & \text{ if } n>0 \text{ is even.} 
\end{cases} 
\end{equation}
which completes the proof.
\end{proof}

\begin{remark}
An alternate proof of the previous theorem can be obtained by computing explicitly
the action of the Lie algebra of global holomorphic vector fields $H^0(F_n;\Theta)$
on $H^1(F_n;\Theta)$ using explicit bases which can be found in \cite[p. 19]{Ca}. 
This approach also permits a better understanding of the normal links 
of the stratifications of $\I$ and $\Jj$, a problem considered in \cite{McD1}. The induced stratification on the link is precisely the cone on the usual stratification of 
the appropriate $\CP^n$ by the secant varieties of the rational normal curve\footnote{The authors 
thank Barbara Fantechi for bringing this stratification of $\CP^n$ to our attention.}.
\end{remark}

\section{The homotopy type of the symplectomorphism group}

In this section we apply the results of the previous sections to study 
 the topology of the symplectomorphism groups of rational
ruled surfaces. 

Let $G_\lambda$ denote the symplectomorphism group of the standard symplectic form
$\omega_\lambda$ on $S^2\times S^2$ or $S^2\tilde \times S^2$, where $\lambda>1$ in
the untwisted case and $\lambda>0$ in the twisted case.

If $G$ is a topological group and $X$ is a $G$-space, we write 
\begin{equation}
\label{hoorbdefn}
X_{hG} = EG \times_G X 
\end{equation}
for the homotopy orbits (or Borel construction) of the $G$-action on $X$.

\subsection{The homotopy decomposition of $G_\lambda$:}

\noindent

\noindent
In this section we write $\Jj_\lambda = \Jj_{\omega_\lambda}$ and $\I_\lambda = \I_{\omega_\lambda}$.
Let $\A_\lambda$ denote the space of almost complex structures tamed by $\omega_\lambda$, $\A_\lambda^i \subset
\A_\lambda$ the subspace of complex structures and $\Omega_{[\lambda]}$ the space of symplectic forms cohomologous to $\omega_\lambda$.

We will need to use the spaces of tame complex structures because the maps between the classifying
spaces $BG_\lambda$ for different values of $\lambda$ are best understood in terms of these spaces
(see \eqref{inclusion} below, cf. \cite[Corollary 2.3]{McD2}). Although maps between the $G_\lambda$
were also defined in \cite[Section 4]{AM}, the construction there
produces only homotopy multiplicative maps and it is not apparent that they induce maps
 on the level of the classifying spaces.

We write $\Jj_{[\lambda]}, \I_{[\lambda]}, \A_{[\lambda]}, \A_{[\lambda]}^i$ for the analogous spaces where the symplectic form $\omega_\lambda$ is replaced with all symplectic forms in the cohomology class of $\omega_\lambda$.

Let 
\[
\K_{[\lambda]} = \{ (J,\omega) \in \Jj_{[\lambda]} \times \Omega_{[\lambda]} \colon J \text{ is compatible with
} \omega \}. \]
$\K_{[\lambda]}^t$ denotes the analogous space with $J$ tamed by $\omega$ and 
$\K_{[\lambda]}^i$ and $\K_{[\lambda]}^{i,t}$ the analogous subspaces where $J$ is required to be integrable.

The above spaces have a natural action of the subgroup of diffeomorphisms
preserving the cohomology class $[\omega_\lambda]$ (cf. \eqref{diffpressymp}). In our situation
this is the subgroup of diffeomorphisms inducing the identity on homology which we denote by
$\Diff_{[0]}$.

Finally, if $S$ is a subset of the isomorphism classes of complex structures, we will decorate
the above spaces with a superscript $S$ to indicate that the corresponding complex structures 
belong to $S$. In the case of not necessarily integrable complex structures the superscript $S$
indicates representability of the appropriate homology class by an embedded $J$-holomorphic curve.

\begin{prop}
\label{equivmoduli}
Let $S$ be a set of isomorphism classes of complex structures. The canonical maps
\[ (\Jj_{\lambda}^S)_{hG_{\lambda}} \to (\Jj_{[\lambda]}^S)_{h\Diff_{[0]}} \]
are weak equivalences. Similarly for tame and/or integrable complex structures. 
\end{prop}
\begin{proof}
We have a $\Diff_{[0]}$ equivariant bundle 
\[ \Jj_\lambda^S \to \K_{[\lambda]}^S \xrightarrow{\pi_2} \Omega_{[\lambda]}.\]
As $\Omega_{[\lambda]} = \Diff_{[0]} /G_\lambda$ by \cite[Theorem 2.4]{LM}, we see that
$E\Diff_{[0]} \times_{\Diff[0]} \K_{[\lambda]}^S = E\Diff_{[0]} \times_{G_\lambda} \Jj_\lambda^S$. The latter space is also a model for $\left(\Jj_\lambda^S\right)_{hG_\lambda}$ hence the inclusion
\[ (\Jj_\lambda^S)_{hG_\lambda} \to (\K_{[\lambda]}^S)_{h\Diff_{[0]}} \]
is a weak equivalence. We have a commutative diagram
\[ 
\xymatrix{
& \K_{[\lambda]}^S \ar[d]^{\pi_2} \\
\Jj_\lambda^S \ar@{^{(}->}[ur] \ar@{_{(}->}[r] & \Jj_{[\lambda]}^S 
}\]
and $\pi_2$ is a $\Diff_{[0]}$-equivariant map with contractible fibers so the result follows
from the homotopy invariance of homotopy orbits.
\end{proof}

\begin{remark}
Using the arguments in the proof of Theorem \ref{Main1} (which corresponds to the case when 
$S$ is the set of all complex structures) one can show that in fact 
all the inclusions between the four different spaces of complex structures associated to a symplectic form restricted to a set $S$ of isomorphism classes induce weak equivalences.
Hence, for each $S$, the eight spaces mentioned in the statement of Proposition \ref{equivmoduli} are all weakly equivalent. We will not need this, however.
\end{remark}

\begin{lemma}
\label{compandtame}
Let $S$ denote a set of isomorphism classes of complex structures on $M$. Then 
\begin{enumerate}[(i)]
\item $\A_{[\lambda]}^{i,S} = \I_{[\lambda]}^S$,
\item $\I_{[\lambda]}^S $ is a finite union of $\Diff_{[0]}$-orbits, each containing 
one of the K\"ahler structures $J_k$ of Theorem \ref{thm:acs} constructed in Section 4.
\end{enumerate}
\end{lemma}
\begin{proof}
Any complex structure $J\in \A_{[\lambda]}^i$ admits a $J$-holomorphic sphere in the homology
class $(1,-k)\in H_2(M;\Z)$ for some $k\in \{0,1,\ldots,\ell\}$ with $\ell<\lambda\leq \ell+1$
as well as a $J$-holomorphic sphere representing the homology class $(0,1)$.
By standard results in complex geometry (see for instance \cite[Proposition V.4.3(i)]{BPV}), it 
follows that $(M,J)$ is complex isomorphic to the Hirzebruch surface $F_{2k}$ or $F_{2k+1}$ 
according to whether we are in the untwisted or twisted cases. Since these
$F_n$ admit a K\"ahler structure with symplectic form in the cohomology class $[\omega_\lambda]$
(constructed in Section 4 above), it follows that $J \in \I_{[\la]}$. This proves (i) and (ii).
\end{proof}

Recall from \cite[Lemma 2.2]{McD2} that $\A_{[\lambda]} \subset \A_{[\lambda+\epsilon]}$ for
$\epsilon>0$. The loop maps\footnote{There is a more geometric construction of homotopy multiplicative maps $G_\lambda \to G_{\lambda+\epsilon}$ in \cite[Section 4]{AM} but it is important for us that these
are actually loop maps, i.e. that they induce maps on the level of classifying spaces. This is not apparent
from \cite[Section 4]{AM} which is why we use the construction in \cite{McD2}.} $G_\lambda \to G_{\lambda+\epsilon}$ are defined by taking the homotopy orbits of $\Diff_{[0]}$ along this inclusion:
\begin{equation}
\label{inclusion}
BG_\lambda \simeq (\A_{[\lambda]})_{h\Diff_{[0]}} \subset (\A_{[\lambda+\epsilon]})_{h\Diff_{[0]}} \simeq BG_{\lambda+\epsilon}. \end{equation}
where the weak equivalences are a consequence of Proposition \ref{equivmoduli} in the case of 
tame complex structures.

Given Theorem \ref{Main1} or, more precisely, Remark \ref{remarktame}, and Proposition \ref{equivmoduli} (in the case of tame integrable complex structures) we may replace $\A_{[\mu]}$ in \eqref{inclusion} by the subspaces $\A_{[\mu]}^i$ of integrable complex structures. By Lemma \ref{compandtame}, this space does
not change when $\mu$ varies in an interval of the form $]k,k+1]$ and so we have the following result of McDuff's.
\begin{prop} \cite[Theorem 1.4]{McD2}
For $k<\la \leq \mu \leq k+1$, the map $BG_\la \to BG_\mu$ is a weak equivalence.
\end{prop}
We can now prove our main result concerning the homotopy type of $BG_\mu$.
First we recall that the \emph{homotopy pushout} of a diagram of spaces 
\[ X \xleftarrow{f} A \xrightarrow{g} Y \]
(or double mapping cylinder of $f$ and $g$) is the quotient space 
\begin{equation}
\label{hopdefn}
  P = \left( X \coprod A \times[0,1] \coprod Y \right)/\sim 
\end{equation}
where $\sim$ denotes the equivalence relation generated by $(a,0)\sim f(a)$ and $(a,1)\sim g(a)$. 
To compute the homology of the homotopy pushout we have the Mayer-Vietoris sequence
\[  \cdots \rightarrow H_*(A) \xrightarrow{(f_*,g_*)} H_*(X) \oplus H_*(Y) \rightarrow  H_*(P) \rightarrow \cdots \]
obtained from the obvious decomposition of $P$ and similarly for cohomology. 
We say that a square
\[ \xymatrix{A \ar[r]^f \ar[d]_g & X \ar[d]^i \\ Y \ar[r]_j & Q } \]
is a \emph{homotopy pushout square} if $if$ and $jg$ are homotopic and the map $P \to Q$ canonically determined
by such a homotopy is a weak equivalence. 

As in the previous section we write $K(n)$ for the connected component of the identity of the K\"ahler isometry group of the Hirzebruch surface $F_n$. Thus
\begin{equation}
\label{defKn}
K(n) = \begin{cases} 
SO(3) \times SO(3) & \text{ if } n=0 \text{ in the untwisted case,} \\
S^1 \times SO(3) & \text{ if } n>0 \text{ is even,} \\
U(2) & \text{ otherwise.}
\end{cases}
\end{equation}
For the sake of simplifying the statement of the following theorem, in the untwisted case, for $0<\lambda\leq 1$ we write $G_\lambda$ for the connected component of the identity of the 
symplectomorphism group $G_1$. By a theorem of Gromov \cite{Gr}, the inclusion $K(0) \subset G_1$ is a weak equivalence and, in the twisted case, $K(1) \subset G_{\lambda}$ is a weak equivalence
for $0<\lambda \leq 1$.

\begin{thm}
\label{maintopological}
Let $\lambda>0$ and let $\ell$ be the integer such that $\ell< \lambda \leq \ell +1$.
There is a homotopy pushout square
\begin{equation}
\label{pushoutsquare}
\xymatrix{
S^{2m-3}_{hK(m)} \ar[d]_j \ar[r]^\pi & BK(m) \ar[d] \\
BG_\lambda \ar[r]_i & BG_{\lambda+1} 
}\end{equation}
where 
\[ m = \begin{cases} 
2\ell + 2 & \text{ in the untwisted case,} \\
2\ell + 3 & \text{ in the twisted case,}
\end{cases}\]
$S^{2m-3}$ is the unit sphere of the representation of $K(m)$ described 
in Theorem \ref{representation}, $\pi$ is the canonical projection and $i$ the map described in \eqref{inclusion}.
\end{thm}
\begin{proof}
Consider the usual stratification of $\I_{\lambda+1}$:
\[ \I_{\lambda+1} = V_0 \cup \ldots \cup V_{\ell+1}. \]
Writing
\[ V_{0k} = V_0 \cup \ldots \cup V_k \]
and $NV_{\ell+1}$ for a tubular neighborhood of $V_{\ell+1}$, we have a homotopy
pushout decomposition
\begin{equation}
\label{pushoutcomplex}
 \xymatrix{ NV_{\ell+1} \setminus V_{\ell+1} \ar[r] \ar[d] & V_{\ell+1} \ar[d] \\
 V_{0\ell} \ar[r] & \I_{\la+1}. }
\end{equation}

Writing $S$ for the set of Hirzebruch surfaces $\{F_i,\ldots,F_{2\ell+i}\}$ where $i=0$ in the untwisted case and $i=1$ in the twisted case, Proposition \ref{equivmoduli} (in the case of 
compatible integrable complex structures) implies that the inclusion 
\[ (V_{0\ell})_{hG_{\la+1}} = (\I_{\la+1}^S)_{hG_{\la+1}} \subset (\I_{[\la+1]}^S)_{h\Diff_{[0]}} \]
is a weak equivalence. By Lemma \ref{compandtame} we have 
\[ \I_{[\la+1]}^S = \A_{[\la+1]}^{i,S} = \A_{[\la]}^i \]
and so another application of Proposition \ref{equivmoduli} to $\I_{\la+1}$ identifies the (homotopy class of the) $G_{\lambda+1}$ homotopy orbits of the bottom row in \eqref{pushoutcomplex} with 
\[ ( \A_{[\la]}^i)_{h\Diff_{[0]}} \subset (\A_{[\la+1]}^i)_{h\Diff_{[0]}} \]
which is exactly the map $i\colon BG_\lambda \to BG_{\la+1}$ in \eqref{inclusion}.

By Theorems \ref{thm:intorbit} and \ref{representation}, the projection $(NV_{\ell+1} \setminus V_{\ell+1}) \to V_{\ell+1}$ has the (weak) homotopy type of the 
projection $G_{\lambda+1} \times_{K(m)} S^{2m-3} \to G_{\lambda+1}/K(m)$ (see Proposition \ref{identtube} for more details).

The only reason why the Theorem doesn't follow immediately by taking homotopy orbits of $G_{\lambda+1}$
on \eqref{pushoutcomplex} is that $NV_{\ell+1}$ is not invariant under the $G_{\lambda+1}$-action.
Nevertheless, the slice theorem implies that for each compact set $W\subset G$ there 
is a smaller tubular neighborhood 
$NV_{\ell+1}'$ sent by the action of $W$ into $NV_{\ell+1}$. This allows us to define
an $A_\infty$-action of $G_{\la+1}$ on $NV_{\ell+1}$ which, as such, is equivalent to the 
left $G_{\la+1}$-action on $G_{\la+1}\times_{K(m)} S^{2m-3}$ (see  
Appendix \ref{quasi} \eqref{nearpushout} for the details). The result follows.
\end{proof}
We note the following immediate consequence of the previous Theorem.
\begin{cor} \cite[Theorem 1.4]{McD2}
\label{connectivity}
With $m$ as in the statement of Theorem \ref{maintopological}, for all $\mu>\lambda$, 
the map $BG_\lambda \to BG_{\mu}$ is $(2m-3)$-connected.
\end{cor}

\subsection{The cohomology groups of $BG_\lambda$:}

\noindent

\noindent
In the rest of the paper we will use Theorem \ref{maintopological} to obtain results concerning
the cohomology of $BG_\lambda$. First we need to compute the effect of the map $\pi$ in \eqref{pushoutsquare} on cohomology.

Observe that $B(S^1\times S^1)=\CP^\infty\times \CP^\infty$ so a choice of basis 
for a $2$-torus $\T$ gives canonical generators $T_1,T_2 \in H^*(B\T;\Z)=\Z[T_1,T_2]$.
Thus the standard bases for the maximal tori of the groups $K(n)$ introduced in Section \ref{deformations} (see \eqref{torusso3}) determine generators for $H^*(BK(n);R)=H^*(B\T;R)^{\Z/2}$  where $\Z/2$ denotes
the Weyl group and $R=\Z[1/2]$ or $\Z$ according to whether $n$ is even or odd. This way we have
\[ 
\begin{cases}
H^*(BK(n);\Z) = \Z[A_n,X_n] & \text{if } n \text{ is odd}, \\ 
H^*(BK(n);\Z[1/2]) = \Z[1/2][A_n,X_n] & \text{ if } n> 0 \text{ is even}
\end{cases}
\]
where $|A_n|=2$ and $|X_n|=4$ while it is easy to see that 
\[ H^*(BK(n);\Z/2) = \Z/2[T,w_2,w_3] \text{ for } n>0 \text{ even,} \] 
with $|T|=|w_2|=2$ and $|w_3|=3$.

\begin{lemma}
\label{eulerclass}
For $n>1$, the Euler class $e_n$ of the isotropy representation of $K(n)$ is given by 
\[ 
e_n = \begin{cases}
\prod_{i=1}^{\frac{n-1}{2}} \left((2i-1)^2 X_n - i(i-1)A_n^2\right) \in H^*(BK(n);\Z) & \text{ if } n
\text{ is odd,} \\
A_n\prod_{i=1}^{\frac{n}{2}-1}(A_n^2-i^2X_n) \in H^*(BK(n);\Z[\frac 1 2]) & \text{ if } n \text{ is even.}
\end{cases}
\]
With any coefficients $R$, the Euler class is a non-zero divisor in $H^*(BK(n);R)$.
\end{lemma}
\begin{proof}
The inclusion of the maximal torus
\[ BS^1\times BS^1 \to BK(n) \]
identifies $H^*(BK(n);R)$ with the ring of polynomials in $H^*(BS^1\times BS^1;R)=R[T_1,T_2]$ invariant
under the Weyl group, where $R=\Z$ or $\Z[\frac 1 2]$ according to whether $n$ is odd or even. 
The Weyl group action switches $T_1$ and $T_2$ for $n$ even and fixes $T_1$ and acts by $-1$ on $T_2$
for $n$ odd. 

If $\chi\colon S^1\times S^1 \to S^1$ is the weight determined by $(k,l)\in \Z^2$, the Euler class
of the line bundle 
\[ E(S^1\times S^1) \times_{S^1\times S^1} \C \to B(S^1\times S^1),\]
where $S^1\times S^1$ acts on $\C$ via $\chi$, is $kT_1+lT_2 \in H^2(B(S^1\times S^1);\Z)$. Since the Euler class is multiplicative, this determines the Euler class of any representation of $K(n)$ once we know the character of the representation on the maximal torus. The formulas for the Euler class now follow from 
Theorem \ref{representation} or, more precisely, from the formulas \eqref{weightsstand} for the
character of the isotropy representations on the maximal tori in terms of the standard weights.

For $n$ odd, $H^*(BK(n);R)$ is a polynomial ring for any $R$ and so it suffices to check that 
$e_n \neq 0$. This is the case because $((2i-1)^2,i(i-1)) = 1$. For $n$ even and any coefficients $R$, $H^*(BK(n);R)$ contains a polynomial ring $R[A_n]$ with $|A_n|$=2 as a retract. Since the coefficient of $A_{n}^{n-1}$ in the formula for $e_{n}$ is $1$, it follows that $e_{n}$ is not a zero divisor.
\end{proof}

\noindent
{\bf Proof of Theorems \ref{cohgps} and \ref{cohgps2}:}
Consider the diagram \eqref{pushoutsquare}. The previous lemma
says that given any coefficient ring R, the map $H^*(\pi;R)$ is the quotient map  
\[ H^*(BK(m);R) \to H^*(BK(m);R)/\langle e_m \rangle. \]
Thus the Mayer-Vietoris sequence associated to \eqref{pushoutsquare} splits into a 
short exact sequence 
\[ 0 \to \langle e_m \rangle = \Sigma^{2m-2} H^*(BK(m);R) \to H^*(BG_{\la+1};R) \to H^*(BG_\la;R) \to 0. \]
When $R$ is a field this splits and since the groups in the extension are finitely generated,
it follows that the same is true over $\Z$.
The statements of Theorems \ref{cohgps} and \ref{cohgps2} follow by induction.

\subsection{Computation of the rational cohomology ring of $BG_\lambda$:}

\noindent

\noindent
Let $\FDiff$ denote the group of fiber preserving diffeomorphisms of 
an $S^2$ bundle over $S^2$ inducing the identity on homology. It will be clear
from the context whether the bundle in question is trivial or not.
To compute the cohomology ring of $BG_\lambda$ we will make use of the loop maps 
\[ G_{\la} \to \FDiff\]
defined in \cite{McD2}.

Writing $\A_{[\infty]}$ for the space of almost complex
structures compatible with some symplectic form, McDuff shows \cite[Proof of Prop 1.1]{McD2}
that the space $\A_{[\infty]}$ is $\Diff_{[0]}$-equivariantly weakly equivalent to the space
$\Diff_{[0]}/\FDiff$ and, hence, taking homotopy orbits we have canonical maps as in \eqref{inclusion}
\[ BG_{\la} \simeq (\A_{[\la]})_{h\Diff_{[0]}} \subset (\A_{[\infty]})_{h\Diff_{[0]}} \simeq B\FDiff. \]
By Theorem \ref{Main1} we can replace $\A_{[\infty]}$ with its subspace $\I_{[\infty]}$ of complex structures.

We begin by noting the following consequence of Theorem \ref{maintopological}.
\begin{prop}
\label{cohomrelation}
For all $\lambda$, the map $BG_{\lambda} \to B\FDiff$ induces a surjection on cohomology
with any coefficients. Moreover 
\[ H^*(B\FDiff) = \lim_{\la} H^*(BG_{\la}). \]
\end{prop}
\begin{proof}
By Lemma \ref{eulerclass} (cf. the proof of Theorems \ref{cohgps} and \ref{cohgps2})
the map $H^*(BG_\mu) \to H^*(BG_\la)$ is surjective for all $\la<\mu$. Since
\[ B\FDiff = \hocolim_{\la} BG_{\la} \]
Lemma \ref{connectivity} implies that the connectivity of the maps $BG_{\mu} \to B\FDiff$ 
tends to $\infty$ with $\mu$. This proves the first statement. The second statement
is obvious.
\end{proof}

The problem that must be overcome in order to compute the cohomology ring of $G_\lambda$
is that of understanding the effect of the map $j$ in \eqref{pushoutsquare} on cohomology 
since, by Lemma \ref{eulerclass},
\begin{equation}
\label{cohomologypushoutsquare}
\xymatrix{ H^\ast S^{2m-3}_{hK(m)}  & \ar[l]_{\pi^\ast} H^\ast BK(m)  \\  H^\ast BG_{\lambda} \ar[u]^{j^\ast} & H^\ast BG_{\lambda+1} \ar[u] \ar[l]^{i^\ast} }
\end{equation}
is a pushout (or equivalently pullback) square of graded abelian groups and hence 
a pullback square of graded rings.

We will do this, for cohomology with rational coefficients, by making use of the 
commutative diagram 
\[
\xymatrix{ S^{2m-3}_{hK(m)}  \ar[r]^\pi \ar[d]_j & BK(m) \ar[d]\ar[rdd] & \\ BG_{\lambda} \ar[rrd] \ar[r]^i & BG_{\lambda+1} \ar[dr] & \\ & & B\FDiff.  }
\]

We begin by analyzing $H^*(B\FDiff;\Q)$. By a Theorem of Smale, the inclusion $SO(3)\subset \Diff^+(S^2)$ is a weak equivalence so we will not distinguish between the two. There is a short exact sequence  
\[ \mathcal{S} \to \FDiff \xrightarrow{e} SO(3) \]
where the map $e$ takes a diffeomorphism to the projection of its action on the base and $\mathcal S$ denotes the gauge group of the appropriate $S^2$ bundle over $S^2$. This yields a fiber sequence
\begin{equation}
\label{fiberseqfdiff}
 B\mathcal S \to B\FDiff \xrightarrow{Be} BSO(3). 
\end{equation}
Now $B\mathcal S$ is the component of the null map in 
\[ \Map(S^2,BSO(3)) \]
in the untwisted case, while in the twisted case it is the component of the essential map
corresponding to the generator of $\pi_2(BSO(3)) = \Z/2$. The latter statement is a consequence
of the fact that, in the twisted case, $\mathcal S$ is the pullback of the diagram of groups
\[ \Map(D^2,SO(3)) \rightarrow \Map(S^1,SO(3)) \leftarrow \Map(D^2,SO(3)) \]
with one of the maps restriction to the boundary and the other restriction followed by conjugation
by a generator of $\pi_1(SO(3))$ together with the fact that the standard simplicial construction of the classifying space functor commutes with fibered products. The two components of $\Map(S^2,BSO(3))$ are equivalent away from the prime $2$ (and hence rationally) since the degree $2$ map of $S^2$ induces an equivalence 
\[ \Map_1(S^2,BSO(3)) \to \Map_0(S^2,BSO(3)). \]

\begin{lemma}
\label{cohfdiffrational}
$H^*(B\FDiff;\Q) = \Q[A,X,Y]$ with $|A|=2$ and $|X|=|Y|=4$, i.e. the rationalization of the 
classifying space $B\FDiff_\Q$ is weakly equivalent to $K(\Q,2)\times K(\Q,4)\times K(\Q,4)$. 
\end{lemma}
\begin{proof}
We have a fiber sequence 
\[ \Omega S^3 \simeq \Omega^2_0 BSO(3) \to \Map_0(S^2,BSO(3)) \to BSO(3).\]
Since $H^*(\Omega S^3;\Q)$ is polynomial generated by a class in degree $2$ and
$H^*(BSO(3);\Q)$ is a polynomial ring generated by a class in degree $4$, the Serre spectral
sequence collapses and so $H^*(B\mathcal S;\Q) = \Q[A,Y]$. Note that, by the discussion preceding 
the statement, this is true in both twisted and untwisted cases. Finally, the spectral sequence of 
\eqref{fiberseqfdiff} also collapses.
\end{proof} 

There is a canonical choice for the degree $2$ generator and one of the degree $4$ generators, namely the pullback of the generator of $H^4(BSO(3);\Q)$ under $Be$ but this is \emph{not} the case for the remaining degree $4$ generator. In order to continue the computation we must choose well defined 
generators in $H^*(B\FDiff;\Q)$. 

By Corollary \ref{connectivity}, in both twisted and untwisted cases, the map 
\[ BG_2 \to B\FDiff \]
is at least $5$-connected and hence induces an isomorphism on cohomology in degrees $\leq 4$.
We will use this fact and the pushout decomposition in Theorem \ref{maintopological} to pick the generators of $H^4(B\FDiff;\Q)$. 

By Theorem \ref{representation}, the groups $K(2)$ and $K(3)$ act transitively on the 
unit sphere of their isotropy representations with isotropy groups $SO(3) \subset K(2)=S^1\times SO(3)$ and $U(1)\times 1 \subset K(3)=U(2)$. Since for $H$ a closed subgroup of $G$ we have $(G/H)_{hG}=EG/H=BH$, in Theorem \ref{maintopological} we have 
\[ S^1_{hK(2)} = BSO(3), \quad S^3_{hK(3)} = BU(1), \]
and hence $BG_2$ is obtained by the homotopy pushouts 
\[ 
\xymatrix{ 
BSO(3) \ar[r]^{\pi} \ar[d]_j & BK(2) \ar[d] & & BS^1 \ar[r]^\pi \ar[d]_j & BK(3) \ar[d] \\
BK(0) \ar[r] & BG_2 & & BK(1) \ar[r] & BG_2}
\]
in the untwisted and twisted cases respectively, with $\pi$ the maps induced by the inclusions
of the isotropy groups. The map $j$ in the left square is the inclusion of the diagonal
by Iglesias' classification of $SO(3)$-equivariant symplectic four-manifolds \cite{I}
(cf. also \cite[Theorem 1.1 (ii)]{AG}) while the map $j$ on the right will be described
below (Proposition \ref{pushouttwisted}). Regardless
of what these maps are, the Mayer-Vietoris sequences of the above diagrams together with 
Corollary \ref{connectivity} imply that the inclusions
\begin{equation}
\label{psimap}
(BSO(3)\times BSO(3)) \vee BS^1 \xrightarrow{\psi} B\FDiff \quad \text{and} \quad 
BK(1)\vee BSU(2) \xrightarrow{\psi} B\FDiff, 
\end{equation}
where $S^1 \subset K(2)=S^1\times SO(3)$ is the inclusion of the first factor and 
$SU(2) \subset K(3)=U(2)$ is the standard inclusion, induce isomorphisms on 
cohomology (even with integral coefficients) in degrees $\leq 4$.

For the rest of this section we will use the following notation for the standard generators 
in $H^*(BK(n);\Q)$:
\begin{equation}
\label{notationgens}
H^*(BK(n);\Q) = \begin{cases} \Q[Y_0,X_0] & \text{ if } n=0,\\
\Q[A_n,X_n] & \text{ if } n>0.
\end{cases} 
\end{equation}
with $|A_n|=2$ and $|X_n|=|Y_0|=4$. Note that $X_0$ corresponds to the second copy of $SO(3)$, the one
which rotates the base of the fibration according to the identification of the groups $K(n)$ in Section 
\ref{deformations}.

\begin{defn}
\label{defgenerators}
Consider the maps $\psi$ in \eqref{psimap}. Let $X,Y \in H^4(B\FDiff;\Q)$ denote the unique classes such that\footnote{We write $X_3$ for the image of the class $X_3 \in H^4(BK(3);\Q)$ in $H^4(BSU(2);\Q)$ and $A_2$ for the image of the class $A_2 \in H^2(BK(2);\Q)$ in $H^2(BS^1;\Q)$ under the obvious inclusions.}
\[ 
\psi^*(X) = \begin{cases} 
X_0 & \text{ in the untwisted case} \\
X_1 & \text{ in the twisted case.}
\end{cases} 
\quad \quad
\psi^*(Y) = \begin{cases} 
Y_0 + A_2^2 & \text{ in the untwisted case}, \\
X_3 & \text{ in the twisted case.}
\end{cases}
\]
and $T \in H^2(B\FDiff;\Q)$ the unique class such that $\psi^*(T)=A_k$, with $k=2$ in the 
untwisted case and $k=1$ in the twisted case.
\end{defn}
\begin{remark}
\label{commentgens}
For $n$ even, $X \in H^4(B\FDiff;\Q)$ is the canonical generator obtained by pulling back 
the generator of $H^4(BSO(3);\Q)$ along the map $Be$ in \eqref{fiberseqfdiff} while $Y$ was chosen to make the formula for $H^*(BG_\lambda;\Q)$ below to agree as much as possible with the one in \cite[Theorem 1.2]{AM}. 

For $n$ odd, since $BK(n) \to B\FDiff \xrightarrow{Be} BSO(3)$ identifies with the projection $BU(2)\to BSO(3)$ and hence sends the generator $Z \in H^4(BSO(3);\Q)$ to $A_n^2-4X_n$, the canonical generator is the class $T^2-4(X+Y) \in H^4(B\FDiff;\Q)$. In this twisted case, this class does not generate $H^4(B\FDiff;\Z)$ and so we decided not to use this as one of the generators for $H^4(B\FDiff;\Q)$.

See also \eqref{otherdefgens} and \eqref{otherdefgens2} below for the relation between the generators in the untwisted case and other geometrically defined classes.
\end{remark}

We now need to compute the effect of the inclusions 
\[ BK(n) \xrightarrow{\psi_n} B\FDiff \]
on rational cohomology. By definition of the generators of $H^*(B\FDiff;\Q)$ 
we know the answer for $n=0,1$ and partly for $n=2,3$. To complete the 
computation we will use the fact that for each $n,m$ with the same parity there are 
$S^1$'s inside $K(n)$ and $K(m)$ which are conjugate inside $\FDiff$. 
This should have an elementary proof but we have only been able to obtain one for 
$(n,m)=(0,2k)$ (which suffices to compute the cohomology ring in the untwisted case).
In order to handle the twisted and untwisted cases uniformly we will take a different tack and use
Karshon's classification of $S^1$-actions to find the conjugate circles inside $G_\lambda$.

We will use the standard bases for the maximal tori $S^1\times S^1 \subset K(n)$ defined in 
Section \ref{deformations} (see \eqref{torusso3}). A circle in $K(n)$ is now described by an integer vector $(a,b) \in \Z^2$. The following result is essentially \cite[Lemma 3]{Ka2} but we include a proof
for the reader's convenience.
\begin{prop}
\label{circlesincommon}
Given $\lambda$ such that there are complex structures compatible with $\omega_\lambda$
isomorphic to $F_k$ and $F_l$, there are $S^1$-equivariant symplectomorphisms between 
\begin{itemize}
\item For $k$ and $l$ odd:
\begin{itemize} 
\item $F_k$ with the $S^1$-action given by $(\frac{l+1}2,\frac{l-1}2)$,
\item $F_l$ with the $S^1$-action given by $(\frac{k+1}2,\frac{k-1}2)$.
\end{itemize}
\item For $k$ and $l$ even:
\begin{itemize}
\item $F_k$ with the $S^1$-action given by $(\frac{l}{2},1)$,
\item $F_l$ with the $S^1$-action given by $(\frac{k}{2},1)$.
\end{itemize}
\end{itemize}
\end{prop}
\begin{proof}
The K\"ahler reduction construction from Section \ref{deformations} provides a standard picture
for the moment polygon of the Hirzebruch surfaces $F_k$ with vertices $(0,0)$, $(1,0)$, $(1,\mu)$ and 
$(0,\mu-n) \in \Lie(\T)^*$ in terms of the moment map basis described in \eqref{momentbasis} for the maximal torus $\T$ of $K(k)$. 

Using the change of basis\footnote{Note that the corresponding change of basis for the tori is the
transpose matrix.} in $\Lie(\T)^*$
\[ \begin{bmatrix} 1 & 0 \\ m & 1 \end{bmatrix} \]
with $m=\pm(k-l)/2$ we can make the slopes of the non-vertical edges of the moment polygons for $k$ and $l$ agree as long as $k$ and $l$ have the same parity. It follows from\cite[Theorem 4.1]{Ka} that $F_k$ and $F_l$ are $S^1$-equivariantly symplectomorphic with respect to the circle actions corresponding to projection onto the $y$-axis (i.e. $(0,1)$) in these bases.

Hence (with $k,l$ of the same parity) the $S^1$'s given in the moment map bases by
\[ \begin{bmatrix} \frac{l-k}2 \\ 1 \end{bmatrix} \text{ on } F_k \quad \quad 
\begin{bmatrix} \frac{k-l}{2} \\ 1 \end{bmatrix} \text {on } F_l \]
produce equivariantly symplectomorphic manifolds (for instance, the two polygons on the left in \cite[Figure 4, p.9]{Ka} correspond to the case $(k,l)=(0,2)$).

Applying the change of basis of Lemma \ref{changeofbasis} to these vectors now completes the proof
(bearing in mind that on $F_0$ the circles written $(a,1)$ and $(-a,1)$ are conjugate).
\end{proof}

\noindent
We can now prove the analog of \cite[Theorem 1.1(ii) and Corollary 4.5]{AG} in the 
twisted case.
\begin{prop}
\label{pushouttwisted}
In the twisted case 
\begin{equation}
\label{decomptwisted}
 BG_2 = \hocolim\left( BU(2) \xleftarrow{B(2,1)} BS^1 \xrightarrow{B(1,0)} BU(2)  \right) 
\end{equation}
Moreover
\[ H^*(BG_2;\Q) = \Q[T,X,Y]/Y(9X-2T^2) \]
and $H^*(BG_2;\Z)$ is the subring generated\footnote{There are infinitely many relations on the
generators of $H^*(BG_2;\Z)$ corresponding to the fact that some elements in $\langle Y(9X-2T^2)\rangle$ are divisible by powers of $3$. These divided classes have to be included as relations so as to not introduce torsion.} over $\Z$ by $T,X,Y$ and $\frac{TY}{3}$.
\end{prop}
\begin{proof}
By Theorem \ref{representation}, the map $BS^1 \xrightarrow{B(1,0)} BU(2)$ is the inclusion of the isotropy group of the representation of $K(3)$ on the normal slice. Proposition \ref{circlesincommon} thus identifies the map $j$ in Theorem \ref{maintopological} with $BS^1 \xrightarrow{B(2,1)} BU(2)$ and the first statement follows.

Using the Mayer-Vietoris sequence of \eqref{decomptwisted} one checks that 
\[
T = (A_1,3A_3) , \quad X = (X_1,2A_3^2), \quad Y = (0,X_3) \in H^*(BK(1);\Z)\times H^*(BK(3);\Z) 
\]
generate $H^*(BG_2;\Z) \subset H^*(BK(1);\Z) \times H^*(BK(3);\Z)$ over $\Q$ and these together with $(0,A_3X_3)$ generate over $\Z$. The result follows.
\end{proof}

\noindent
We will need the following simple computations which are left as an exercise. 
\begin{lemma}
\label{simplecomputations}
Writing $H^*(BS^1)=\Q[T]$, $H^*(BS^1\times BSO(3);\Q) = H^*(BU(2);\Q) = \Q[A,X]$ and 
$H^*(BSO(3)\times BSO(3);\Q) = \Q[X,Y]$, the map $S^1 \xrightarrow{(a,b)} S^1\times SO(3)$ induces the map 
\[ A \mapsto aT, \quad \quad \quad X \mapsto b^2 T^2, \]
$S^1 \xrightarrow{(a,b)} U(2)$ induces  
\[ A \mapsto (a+b)T, \quad \quad  \quad X \mapsto (ab)T^2.\]
and $S^1 \xrightarrow{(a,b)} SO(3)\times SO(3)$ induces
\[ X \mapsto a^2T^2, \quad \quad \quad Y \mapsto b^2T^2. \]
\end{lemma}

\noindent
We are now in a position to understand the effect on cohomology of the inclusions of the groups
$K(n)$ in $\FDiff$ which will be the crucial input for the computation of the rational cohomology rings.
\begin{prop}
\label{mapsoncohomology}
Consider the inclusions $BK(n)\xrightarrow{\psi_n} B\FDiff$. 
\begin{itemize}
\item If $n>0$ is even then
\begin{eqnarray*}
\psi_n^*(T) &  =  & \frac{n}{2} A_n \\
\psi_n^*(X) & =  & X_n \\
\psi_n^*(Y) & =  & A_n^2 + \frac{n^2}{4} X_n. \\
\end{eqnarray*}
\item If $n>0$ is odd then
\begin{eqnarray*}
\psi_n^*(T) & = & n A_n \\
\psi_n^*(X) & = & \frac{n^2-1}{4} A_n^2 + (1-\frac{n^2-1}{8}) X_n \\
\psi_n^*(Y) & = & \frac{n^2-1}{8} X_n.
\end{eqnarray*}
\end{itemize}
\end{prop}
\begin{proof}
For $n=2$, we only need to compute the coefficients of $X_2$ in $\psi_2^*(X)$ and $\psi_2^*(Y)$
(since the coefficients of $A_2^2$  are determined by definition of the generators $X$ and $Y$).
By Proposition \ref{circlesincommon}, the groups $K(0)$ and $K(2)$ contain a circle in common
written in the standard bases as $(1,1)$ and $(0,1)$. Therefore we have a homotopy commutative diagram
\[ \xymatrix{ BS^1 \ar[r]^{B(0,1)} \ar[d]_{B(1,1)} &  BK(2) \ar[d]^{\psi_2} \\
 BK(0) \ar[r]_{\psi_0} & B\FDiff }\]
and so Lemma \ref{simplecomputations} implies that these coefficients are both $1$
as the statement indicates.

For $n=3$, we only need to find the coefficient of $A_3$ in $\psi_3^*(T)$ and the 
coefficients of $A_3^2$ in $\psi_3^*(X)$ and $\psi_3^*(Y)$.
By Proposition \ref{circlesincommon} we have a commutative diagram
\begin{equation}
\label{pushoutwisted}
\xymatrix{ BS^1 \ar[r]^{B(1,0)} \ar[d]_{B(2,1)} &  BK(3) \ar[d]^{\psi_3} \\
 BK(1) \ar[r]_{\psi_1} & B\FDiff  }
\end{equation}
and so Lemma \ref{simplecomputations} implies that $\psi_3^*(T)=3A_3$, $\psi_3^*(X)=2A_3^2$ and $\psi_3^*(Y)=X_3$ as the statement indicates.

For all higher $n$, Proposition \ref{circlesincommon} and the previous arguments tell us the
effect on cohomology of $BS^1 \xrightarrow{(a,b)} BK(n) \xrightarrow{\psi_n} B\FDiff$ for two independent vectors $(a,b)$ and so simple algebra yields the remaining formulas.
\end{proof}

\begin{cor}
\label{kernel}
The kernel of the map $\psi_n^*: H^*(B\FDiff;\Q) \to H^*(BK(n);\Q)$ is the ideal
\begin{enumerate}[(i)]
\item $\langle T\rangle$ if $n=0$,
\item $\langle \frac{n^4}{16} X - \frac{n^2}{4}Y+ T^2\rangle$ if $n$ is even,
\item $\langle \frac{(n^2-1)n^2}{8}(X+Y) -n^2 Y - \frac{(n^2-1)^2}{32}T^2 \rangle$ if $n$ is odd.
\end{enumerate}
\end{cor}

\noindent

The following result collects the statements concerning the rational cohomology
ring in Theorems \ref{Main2} and \ref{Main3} (see Remark \ref{concrem} concerning the formulas in the
twisted case):
\begin{thm} 
\label{Main23}
Let $\ell <\lambda \leq \ell+1$. With the choice of generators indicated in Definition 
\ref{defgenerators}, the map $H^*(B\FDiff;\Q) \to H^*(BG_{\la};\Q)$ is the quotient
map
\[ \Q[T,X,Y] \longrightarrow \Q[T,X,Y]/(R_{\ell}(T,X,Y)) \]
where
\[ R_{k}(T,X,Y)= \begin{cases}
T(X-Y+T^2)\ldots(k^4X-k^2Y+T^2) & \text{ in the untwisted case,}\\
Y\ldots\left((2k+1)^2(\frac{k(k+1)}2 U -Y)-\frac{k^2(k+1)^2}{2}T^2\right)  &\text{ in the twisted case.}
\end{cases}
\]
and $U=X+Y$.
\end{thm}
\begin{proof}
The proof is by induction. The result is clear for $\ell=0$. Assume $\ell \geq 0$ and the result holds for $\lambda\in]\ell,\ell+1]$. Then by Theorem \ref{maintopological} and  Lemma \ref{eulerclass} 
we have a pullback diagram of rings 
\begin{equation}
\label{mainsquare}
\xymatrix{ \Q[A_m,X_m]/\langle e_m \rangle  & \ar[l]_{\pi^\ast} \Q[A_m , X_m] \\
\Q[T,X,Y]/\langle R_{\ell}(T,X,Y)\rangle \ar[u]^{j^\ast} & H^\ast(BG_{\lambda+1};\Q) = \Q[T,X,Y]/I_{\ell+1} \ar[u] \ar[l] }
\end{equation}
where we have used Proposition \ref{cohomrelation} to express $H^*(BG_{\lambda+1};\Q)$ 
as a quotient of $H^*(B\FDiff;\Q)$ by an ideal $I_{\ell+1}$, $e_m$ denotes the 
Euler class calculated in Lemma \ref{eulerclass} and $m=2\ell+2$ or $2\ell+3$ according to whether we are in the untwisted or twisted case.

We must have $I_{\ell+1} \subset \langle R_{\ell}  \rangle$. On the other hand, 
by Corollary \ref{kernel}, we also have
\begin{equation}
\label{secondcontainment}
I_{\ell+1} \subset \begin{cases}
\langle \frac{m^4}{16} X - \frac{m^2}{4}Y+ T^2\rangle  & \text{ if } m \text{ is even,}\\
\langle \frac{(m^2-1)m^2}{8}(X+Y) -m^2 Y - \frac{(m^2-1)^2}{32}T^2 \rangle & 
\text{ if } m \text{ is odd.}
\end{cases}
\end{equation}

Since $R_{\ell}$ and the polynomials $k_m(T,X,Y)$ appearing in \eqref{secondcontainment} are coprime it follows that 
\[ I_{\ell+1} \subset \langle R_{\ell} k_m \rangle = \langle R_{\ell+1} \rangle. \]
Denoting by $d$ the degree of $R_\ell$, the pullback square \eqref{mainsquare} gives the following generating function for the graded ring $H^*(BG_{\lambda+1};\Q)$:
\[ \chi = \frac{1}{(1-t^2)(1-t^4)} + \frac{1-t^d}{(1-t^2)(1-t^4)^2} - \frac{1-t^d}{(1-t^2)(1-t^4)}. \]
This simplifies to 
\[ \chi = \frac{1-t^{d+4}}{(1-t^2)(1-t^4)^2}, \]
which is the generating function for $\Q[T,X,Y]/\langle R_{\ell+1}(T,X,Y)\rangle$. Hence
$I_{\ell+1}= \langle R_{\ell+1} \rangle$
\footnote{Checking that $\psi_m^*$ does indeed send $R_\ell$ to the ideal generated by $e_m$ (so
that $j^*$ is well defined) is a recommended confidence building activity.}.
\end{proof}

\begin{remark} 
\label{concrem}
Making the change of variables
\[  z = T, \quad x = 4U - T^2, \quad y = 4U + 32Y - 2T^2 \]
we have 
\[ (2k+1)^2\left(\frac{k(k+1)}2 U -Y\right)-\frac{k^2(k+1)^2}{2}T^2 = \frac{1}{32}\left(-z^2+ (2k+1)^4x - (2k+1)^2y\right) \]
and so using the generators $x,y,z$ we obtain the following formula for $H^*(BG_\lambda;\Q)$ in the statement of Theorem \ref{Main3} which is similar to the formula in the untwisted case:
\[ 
H^*(BG_\lambda,\Q) = \frac{\Q[x,y,z]}{ \langle \prod_{i=0}^\ell(-z^2+(2i+1)^4 x - (2i+1)^2 y) \rangle} 
\quad \text{ for } 0\leq \ell< \lambda \leq \ell+1.
\]
The generators $x$ and $z$ have natural geometric descriptions (see Remark \ref{commentgens}) but we have no geometric interpretation for $y$.
\end{remark}

\begin{remark}
\label{Dusamistake}
The ring structure obtained in Theorem \ref{Main23} and Remark \ref{concrem}
differs from the one previously calculated in \cite[Theorem 1.2 and Theorem 1.5]{AM}. 
There are two different reasons for the difference. 

Regarding the multiples of $T^2$ in the factors that make up the 
relation, the problem can be traced to a misapplication of \cite[Theorem 5.4]{AA} in
\cite[p. 1007]{AM}. Using the notation of \cite{AA}, $\overline K(d\mu)$ depends
only on the value of $d\mu$ in the associated graded vector space of the filtration
of the Sullivan model by word length. Thus the higher Whitehead products
provide information only on the image of the relation in this associated graded vector
space and have no bearing on the coefficients of the decomposable $T^2$ terms 
(which have higher filtration).

In this way one sees that the higher Whitehead products 
in $\pi_*(G_\lambda)\otimes \Q$ can not be used exclusively to compute the ring structure in $H^*(BG_\lambda;\Q)$.

The explanation for the remaining difference in the twisted case (regarding the 
coefficients of $X$ and $Y$ in the factors that make up the relation) lies in 
a mistake in \cite[Lemma 2.11]{AM}. Indeed, it follows easily from 
Proposition \ref{mapsoncohomology} that (with the notation of \cite{AM}) we 
have 
\begin{eqnarray*}
\alpha_k & = & (2k+1) \alpha_0 \in \pi_1(G^1_\lambda),\\
\xi_k & = & \xi_0 + \frac{k(k+1)}{2} \eta \in H_3(G^1_\lambda;\Z).
\end{eqnarray*}

\end{remark}

\subsection{Computation of the cohomology ring away from 2 in the untwisted case:}

\noindent

\noindent
In this section, we calculate $H^*(BG_{\lambda};\Z[1/2])$ in the untwisted case. 
This will be done by combining Proposition \ref{cohomrelation} with the calculation of
$H^*(B\FDiff;\Z[1/2])$. Henceforth, all spaces will be localized away from the prime $2$.
We will write $R=\Z[1/2]$.

Since we are working away from the prime $2$, it is easy to see that, in the untwisted case, 
$\FDiff$ is equivalent to the semi-direct product
\[ \mathcal{G} = SO(3) \ltimes \Map(S^2,SU(2)) \]
where $SO(3)$ acts on $\Map(S^2,SU(2))$ by pre-composition of its standard action on $S^2 = \C P^1$. 
\begin{thm} \label{main}
Let $R=\Z[1/2]$. $H^*(B\mathcal G;R)$ is a free module over $R[x,y]$ on generators $a_k,b_k$
with $k \geq 0$:
\[ H^*(B\mathcal{G},R) = R[x,y]\langle a_0,b_0,a_1,b_1,a_2,\ldots \rangle, \]
where $a_0 = 1$, $|x|=|y|=4$, $|b_k|=4k+2$, and $|a_k|=4k$. Moreover, $H^*(B\mathcal{G},R)$ is isomorphic to the subring of $\Q[x,y,z]$, with $|z|=2$, when $b_k$ and $a_k$ are identified respectively with:
\[ \frac{z}{(2k+1)!}\prod_{i=1}^k(z^2+ i^4 x-i^2 y), \quad \frac{z^2}{(2k)!}\prod_{i=1}^{k-1}(z^2+ i^4 x-i^2 y).
\]
\end{thm}

\begin{remark}
One can see that the groups of fiber preserving diffeomorphisms for the twisted and untwisted bundles are equivalent away from the prime $2$, and so the previous Theorem describes 
$H^*(B\FDiff;\Z[1/2])$ also in the twisted case. We will not use this, however. 
\end{remark}

\subsection{Proof of Theorem \ref{main}:}

\noindent

\noindent
Recall that $\mathcal{G} = SO(3) \ltimes \Map(S^2,SU(2))$. $\mathcal{G}$ contains a subgroup $G = SO(3) \times SU(2)$ extending the group of constant maps. Let $\hT \times S^1$ be the maximal torus of $G$. Notice that $\hT$ acts on $\Map(S^2,SU(2))$ by pre-composition with the action of $\hT$ on $S^2$ given by rotation about the vertical axis. $S^1$ is seen as the subgroup of constant maps with value in the maximal torus of $SU(2)$. Let $ \tau_1$ and $\tau_2$ be elements in each factor of $G$  that map to generators of the Weyl group. Notice also that $\mathcal{G}$ contains the $\hT$ invariant subgroup $\Omega^2 SU(2) \subset \Map(S^2,SU(2))$ consisting of maps that take the north pole to the identity element. Here and henceforth, we will fix the north pole of $S^2$ as the basepoint.

\bigskip

\noindent
The proof of Theorem \ref{main} uses a sequence of inclusions of subgroups:
\[ \Omega^2 SU(2) \subset \mathcal{K} \subset \mathcal{H} \subset \mathcal{G} \]
\noindent 
that induces maps of classifying spaces:
\[ \Omega SU(2) \longrightarrow B\mathcal{K} \longrightarrow B\mathcal{H} \longrightarrow B\mathcal{G} \]
with $B\mathcal{H}$ equivalent to the homotopy orbit space of a $\Z/2$-action on $B\mathcal{K}$ and $B\mathcal{G}$ equivalent to the homotopy orbit space of a $\Z/2$-action on $B\mathcal{H}$. More precisely, we have:

\begin{defn}\label{G1}
$\mathcal{H} \subset \mathcal{G}$ is defined as the subgroup 
\[ \mathcal{H} = \hT \ltimes \Map(S^2,SU(2)) \subset SO(3) \ltimes \Map(S^2,SU(2)) = \mathcal{G} \]
We define an involution on $\mathcal{H}$ induced by conjugation with the element $\tau_1 \in {\mathcal G}$. This action is given by inversion on the $\hT$ factor, and by the action induced on $\Map(S^2,SU(2))$ via the action of $\tau_1 \in SO(3)$ by left multiplication on $S^2 = SO(3)/\hT$. This induces an involution on $B\mathcal{H}$, which we also denote by $\tau_1$.
\end{defn}
\begin{defn}\label{G2}
$\mathcal{K} \subset \mathcal{H}$ is defined as the subgroup 
\[ \mathcal{K} = (\hT \times S^1) \; \ltimes \; \Omega^2 SU(2) \subset (\hT \times SU(2)) \; \ltimes \; \Omega^2 SU(2) = \hT \ltimes \Map(S^2,SU(2))  = \mathcal{H} \]
where $S^1 \ltimes \Omega^2 SU(2)$ may be seen as the subspace of maps from $S^2$ to $SU(2)$ that map the basepoint of $S^2$ to $S^1$. Define an involution on $\mathcal{K}$ induced by conjugation with the element  $\tau_2 \in \mathcal{H}$. This action preserves the $\hT$ factor, acts by inversion on the $S^1$-factor, and acts by pointwise conjugation with $\tau_2 \in SU(2)$ on $\Omega^2SU(2)$. As before, this induces an involution on $B\mathcal{K}$ denoted by $\tau_2$. 
\end{defn}
\noindent
>From the above descriptions, we can describe the homotopy type of the respective classifying spaces away from the prime $2$:
\begin{align}
& B\mathcal{K} = E(\hT \times S^1) \times_{\hT \times S^1} \Omega^2 BSU(2)  \label{co1} \\
& B\mathcal{H} = E \hT \times_{\hT} \Map(S^2,BSU(2)) \label{co2}  \\
& B\mathcal{G} = E SO(3) \times_{SO(3)} \Map(S^2,BSU(2))  \label{co3}
\end{align}
\noindent
It is a standard argument to identify the invariant cohomology rings
\[ H^*(B\mathcal{K},R)^{\tau_2} = H^*(B\mathcal{H},R), \quad \mbox{and} \quad H^*(B\mathcal{H},R)^{\tau_1} = H^*(B\mathcal{G},R). \]

\bigskip

\noindent
The action of $\tau_1$ is subtle. To understand this action, we start by $\hT$- equivariantly decomposing $S^2$ as a pushout of two hemispheres intersecting over the equator. From \eqref{co2} we get a pullback diagram:
\[
\xymatrix{ B\mathcal{H}  \ar[r]^{ev_N} \ar[d]^{ev_S} & B\hT \times BSU(2) \ar[d] \\  B\hT \times BSU(2) \ar[r] & E\hT \times_{\hT} LBSU(2) }
\]
where $ev_N$ and $ev_S$ denote evaluation at the north and south pole. Also, the notation $LBSU(2)$ refers to the free loop space of $BSU(2)$ i.e. $\Map(S^1,BSU(2))$. 
We note that the $\tau_1$-action on $B\mathcal{H}$ described above has the property of  swapping the corners of the pullback diagram since it interchanges the north and south pole of $S^2$, and inducing an action on the equator given by inversion. 

\bigskip

\noindent
Notice that \eqref{co1} shows that the cohomology of $B\mathcal{K}$ is equivalent to the equivariant cohomology of $\Omega SU(2)$. This calculation has been made in  \cite{HHH} (Section 6.2). 

\medskip
\noindent
Before we give a description of this cohomology, let us set some notation. Let $u,v \in H^2(B\hT \times BS^1,\Z)$ be the canonical generators corresponding to the two factors respectively. It also follows from an easy spectral sequence argument that $H^2(B\mathcal{G},\Z)$ is a free $\Z$ module generated by a unique class $w$ that restricts to the generator of $H^2(\Omega SU(2), \Z)$. Moreover, $w$ restricts trivially to the cohomology of $B(\hT \times S^1)$ since the inclusion of $\hT \times S^1 \subset \mathcal{G}$ factors through the group $SU(2) \times SU(2)$. Hence, $w$ restricts to the canonical generator of $\hT \times S^1$-equivariant cohomology of $\Omega SU(2)$ described in \cite{HHH}. Therefore, we have:

\begin{thm}\cite[Section 6.2]{HHH}
Let $u,v$ and $w$ be the classes defined above. Then, the cohomology of $B\mathcal{K}$ with coefficients in the ring $\Z$ is given by the the following free module over the ring $\Z[u,v]$ on generators $f_k$, $g_k$, $k \geq 0$:
\[ H^*(B\mathcal{K},\Z) = \Z[u,v]\langle g_0,f_0,g_1,f_1,g_2,\ldots \rangle, \quad \mbox{where} \quad g_0 = 1 \]
where the degree of $f_k$ is $4k+2$, and that of $g_k$ is $4k$. Moreover, as a ring, we may identify $H^*(B\mathcal{K},\Z)$ as the subring of $\Q[u,v,w]$, where the degree of $Z$ is $2$, and the elements $f_k$ and $g_k$ are identified respectively to the elements:
\[ \frac{w}{(2k+1)!}\prod_{i=1}^k((w+ i^2u)^2 -4i^2 v^2), \quad \frac{w(w + k^2u+2kv)}{(2k)!}\prod_{i=1}^{k-1}((w+ i^2u)^2 -4i^2 v^2).
\]
\end{thm}

\bigskip

\noindent
>From the previous remark, and the description in \cite{HHH} we see that the $\tau_2$ action has the property:
\[  \tau_2(u) = u, \quad \tau_2(v) = -v, \quad \tau_2(w) = w \]
Filtering $H^*(B\mathcal{K},R)$ by powers of $v$, and taking $\Z/2$ invariants, we easily derive the following result
\begin{prop}
The cohomology of $B\mathcal{H}$ with coefficients in the ring $R$ is given by the the following free module over the ring $R[u,v^2]$ on generators $b_k$, $a_k$, $k \geq 0$:
\[ H^*(B\mathcal{H},R) = R[u,v^2]\langle a_0,b_0,a_1,b_1,a_2,\ldots \rangle, \quad \mbox{where} \quad a_0 = 1 \]
where the degree of $b_k$ is $4k+2$, and that of $a_k$ is $4k$. Moreover, as a ring, we may identify $H^*(B\mathcal{H},R)$ as the subring of $\Q[u,v,w]$, where the degree of $w$ is $2$, and the elements $b_k$ and $a_k$ are identified respectively to the elements:
\[ \frac{w}{(2k+1)!}\prod_{i=1}^k((w+ i^2u)^2 -4i^2 v^2), \quad \frac{w^2}{(2k)!}\prod_{i=1}^{k-1}((w+ i^2u)^2 -4i^2 v^2).
\]
\end{prop}
\noindent
The hard part now is to identify the action of $\tau_1$ on $H^*(B\mathcal{H},R)$.
\begin{prop}\label{action}
The action of $\tau_1$ on $H^*(B\mathcal{H},R)$ is given by 
\[ \tau_1(w) = w, \quad \tau_1(u) = -u, \quad \tau_1(v^2) = v^2 - uw \]
\end{prop}
\noindent
Before we proceed with the proof of the Proposition \ref{action}, let us see how we may derive Theorem \ref{main} for the cohomology of $B\mathcal{G}$ from this action. Observe that the following elements are invariant under $\tau_1$:
\[ w, \quad u^2, \quad 2v^2-uw \]
Notice also that we have the following equality:
\begin{equation}
\label{formulamod}
(w+i^2u)^2-4i^2v^2 = w^2 + i^4u^2 -i^22(2v^2-uw) 
\end{equation}
It follows that all the elements $a_k, b_k$ are invariant under $\tau_1$. We claim:
\begin{prop}
\[ H^*(B\mathcal{G},R) = H^*(B\mathcal{H},R)^{\Z/2} = R[u^2,2v^2-uw]\langle a_0,b_0,a_1,b_1,\ldots \rangle \]
\end{prop}
\begin{proof}
Since $u,w$ are elements of $H^*(B\mathcal{H},R)$, we may replace the element $v^2$ by $2v^2-uw$ to get:
\[ H^*(B\mathcal{H},R) = R[u,2v^2-uw]\langle a_0,b_0,a_1,b_1,\ldots \rangle \]
\noindent
The proof follows on filtering $H^*(B\mathcal{H},R)$ by powers of $a$, and taking invariants. 
\end{proof}
\noindent
Over $R$ we can replace the generator $2v^2-uw$ with $2(2v^2-uw)$ hence, taking note of 
\eqref{formulamod} and setting 
\begin{equation}
\label{otherdefgens}
x = u^2, \quad y=2(2v^2-uw), \quad z=w
\end{equation}
we obtain
\[ H^*(B\mathcal{G},R) =  R[x,y]\langle a_0,b_0,a_1,b_1,a_2,\ldots \rangle, \quad \mbox{where} \quad 
a_0 = 1. \]
This completes the proof of Theorem \ref{main} assuming Proposition \ref{action}.
\begin{proof}[Proof of Proposition \ref{action}] 
It is clear from the definition that the action must preserve $w$. It also follows from Definition \ref{G1} that the action reverses the sign of $u$. Hence, only the action on $v^2$ needs to be described. Notice that $v^2 = ev_N^*(\sigma)$, where $\sigma \in H^4(BSU(2),\Z)$ is a generator. Recall that the $\Z/2\langle \tau_1 \rangle$ action on the diagram defining $B\mathcal{H}$ as a pullback has the property of switching the corners. Hence, it follows that $\tau_1(v^2) = ev_S^*(\sigma)$. We reconsider the pullback:
\[
\xymatrix{ B\mathcal{H}  \ar[r]^{ev_N} \ar[d]^{ev_S} & B\hT \times BSU(2) \ar[d] \\  B\hT \times BSU(2) \ar[r] & E\hT \times_{\hT} LBSU(2) }
\]
and consider the Serre spectral sequence for the right vertical map, seen as a fibration, whose fiber is $\Omega SU(2)$. Let $\alpha$ be the element in the $E_2$-term representing a generator of $H^2(\Omega SU(2),\Z)$. It is well known \cite{K,ABKS} that the integral cohomology of $E\hT\times_{\hT} LBSU(2)$ is free in degree $2$, generated by $u$, and trivial in degree $5$. Hence the class $\pm \alpha u$ represents $\sigma$ in the cohomology of the total space $B\hT \times BSU(2)$. Recall that the class $w \in H^2(B{\mathcal H},R)$ is represented by $\alpha$ in the Serre spectral sequence for the left vertical fibration (up to an indeterminacy given by a multiple of $u$). Hence, by choosing a suitable sign for $w$, we notice that $\tau_1(v^2) = -uw$ modulo lower filtration in the Serre spectral sequence of the left vertical map, seen as a fibration. We may therefore write:
\[ \tau_1(v^2) = -uw + av^2 + bu^2 \]
applying $\tau_1$ again to this equation tells us that $a=1,b=0$.
\end{proof}

\noindent
We can now prove Theorem \ref{Main2}.

\subsection{Proof of Theorem \ref{Main2}:}

\noindent 

\noindent
The canonical projection map 
\[ B\phi: B\mathcal{G} \longrightarrow B\FDiff \]
is an isomorphism on cohomology with $\Z[1/2]$ coefficients. Comparing the definition 
of the classes $x,y$ and $z$ in the rational cohomology of $B\mathcal{G}$ given in \eqref{otherdefgens} above, with that of $X,Y$ and $T$ in Definition \ref{defgenerators}, we see that 
\begin{equation}
\label{otherdefgens2}
B\phi^*(X)=x, \quad B\phi^*(Y)=y, \quad  B\phi^*(T)=z. 
\end{equation}
Henceforth we identify these graded
rings by the map $B\phi^*$.

Recall that we have an inclusion map 
\[  BK(2k) \xrightarrow{\psi_{2k}} B\FDiff \] 
whose effect on rational cohomology was described in Proposition \ref{mapsoncohomology}.
We recall that
\[ \psi_{2k}^*(X) = X_{2k}, \quad \psi_{2k}^*(Y) = k^2 X_{2k} + A_{2k}^2,
 \quad \psi_{2k}^*(T) = kA_{2k} \]
where we are using the notation established in \eqref{notationgens}. 
Consider the classes $a_k, b_{k-1}$ in the cohomology of $B\FDiff$ defined in Theorem \ref{main}. Then
\[ \psi_{2k}^*(b_{k-1}) = A_{2k} \prod_{i=1}^{k-1}(A_{2k}^2-i^2 X_{2k}) = e_{2k}, \quad \quad \psi_{2k}^*(a_k) = \frac{A_{2k} \, e_{2k}}{2} \]
where $e_{2k}$ denotes the Euler class calculated in Lemma \ref{eulerclass}. Hence, $B\psi_{2k}^*$ maps the submodule $R[x,y]\langle b_{k-1}, a_k \rangle$ isomorphically onto the ideal generated by the Euler class $e_{2k}$. Moreover, it is also clear that the classes $b_i$ and $a_j$ map to zero if $i\geq k$ and $j > k$. It now follows by induction using Theorem
\ref{maintopological} and Proposition \ref{cohomrelation} that  
the kernel of the map 
\[ H^*(B\FDiff,R) \rightarrow H^*(BG_l,R)\] 
is the submodule generated over $R[x,y]$ by the elements $b_i, a_j$ where $i\geq l$ and $j>l$. Furthermore, one has the following identification:
\[ H^*(BG_\lambda,\Q) =  \Q[x,y]\langle a_0,b_0,a_1,\ldots a_l \rangle = \frac{\Q[x,y,z]}{\langle z\prod_{i=1}^l(z^2+i^4 x - i^2 y)\rangle } \]
and so we may identify $H^*(BG_\lambda,R) = R[x,y]\langle a_0,b_0,a_1,\ldots a_l \rangle$ naturally as a subring of the above quotient. This completes the proof of Theorem \ref{Main2}.

\newpage

\appendix

\section{$\db$-operators and the derivative of the Nijenhuis tensor}
\label{app:A}
\noindent
In this appendix we recall some standard almost-complex geometry facts and 
prove a formula for a certain derivative of the Nijenhuis tensor (Corollary~\ref{cor:dN02}). This formula, needed in Section $2$, is elementary 
and probably well-known to experts, but we were unable to find it in the literature.

Let $(M,J)$ be an almost-complex manifold. As usual, we identify $TM$ with the $+i$ eigenspace of the action of $J$ on $TM\otimes \C$ and consider the decomposition of
$\Omega(M)$, the space of complex valued differential forms on $M$, according to
$(p,q)$-type. In particular, $\Om^1(M) = \Om^{1,0}_J (M) \oplus \Om^{0,1}_J (M)$
and $\Om^2(M) = \Om^{2,0}_J (M) \oplus \Om^{0,2}_J (M)\oplus \Om^{1,1}_J (M)$. 

For a complex valued function $f\in\Om^0 (M)$ one defines 
$\db f \in \Om^{0,1}_J (M)$ as 
\[
\db f = (df)^{0,1}\,.
\]
Similarly, given $\al\in\Om^{0,1}_J (M)$, one defines $\db\al\in\Om^{0,2}_J (M)$ as
\[
\db \al = (d\al)^{0,2}\,.
\]
One easily checks that this $\db$-operator, 
$\db:\Om^{0,1}_J (M)\to \Om^{0,2}_J (M)$, satisfies the Leibnitz rule
\[
\db(f\al) = (\db f)\wedge\al + f\db\al\,,\ \forall f\in\Om^0
(M)\,,\ \al\in\Om^{0,1}_J (M)\,.
\]

One can also define an appropriate $\db$-operator on $\Omega(TM)$, the space of $TM$-valued differential forms on $M$. This operator has a particularly simple
expression on $\Omega^0 (TM)$ involving the Nijenhuis tensor (cf.~\cite[(2.6.3),(2.7.1)]{Ga}).

\begin{defn} \label{defn:nijenhuis}
Given an almost-complex manifold $(M,J)$, the Nijenhuis tensor $N_J$
is defined by
\[
N_J (X,Y) = [JX,JY] - J[X,JY] - J[JX,Y] - [X,Y] \,,\ \forall
X,Y\in\Om^0 (TM)\,.
\]
Note that $N_J \in \Om_J^{0,2} (TM) = \Om_J^{0,2} (M)\otimes
\Om^{0}(TM)$, i.e.
\[
N_J(JX,Y) = N_J(X,JY) = - J N_J (X,Y)\,,\ \forall
X,Y\in\Om^0 (TM)\,.
\]
\end{defn}

\begin{defn} \label{defn:db0TM}
Define the operator
$\db:\Om^0(TM) \rightarrow \Om^{0,1}(TM)$, $Y \mapsto \db Y$,
by
\[
(\db Y)(X) \equiv \db_X Y \equiv \frac{1}{2} \left\{
[X,Y] + J[JX,Y] + \frac{1}{2} N_J(X,Y)\right\}\,,\ \forall
X,Y\in\Om^0(TM) \,. 
\]
One easily checks that $\db Y \in \Om^{0,1}(TM)$ is indeed well-defined 
by this formula, i.e. $(\db Y)(X)$ is a tensor in $X$ and 
$(\db Y)(JX) = -J (\db Y)(X)$.
\end{defn}

The following proposition shows that important properties of an integrable 
$\db$-operator are still valid in this non-integrable context.

\begin{prop} \label{prop:db0TM}
The operator $\db:\Om^0(TM) \rightarrow \Om^{0,1}(TM)$ has the following
properties:
\begin{itemize}
\item[(i)] $\db (f\cdot Y) = (\db f)\otimes Y + f \cdot (\db Y)$,
  for any function $f\in\Om^0(M)$. 
\item[(ii)] $\db_X(JY) = J \db_X Y$.
\end{itemize}
\end{prop}
\begin{proof}
Property (i) follows from the following calculation:
\begin{align}
2 \db (fY) (X) & = [X,fY] + J[JX,fY] + \frac{1}{2} N_J(X,fY) \notag \\
& = (X\cdot f)Y + J[(JX)\cdot f]Y + 2f \db_X (Y) \notag \\
& = [(df)(X) + i(df)(JX)]Y + 2f \db_X (Y) \notag \\
& = [2(\db f)(X)] Y + 2f \db_X (Y)\,.\notag
\end{align}
To prove property (ii) note that
\begin{align}
2[\db_X Y + J\db_X (JY)] & = [X,Y] + J[JX,Y] + \frac{1}{2} N_J(X,Y)
\notag \\
& \quad + J[X,JY] - [JX,JY] + \frac{1}{2} J N_J (X,JY) \notag \\
& = [X,Y] + J[JX,Y] + J[X,JY] - [JX,JY] \notag \\
& \quad + \frac{1}{2} (N_J(X,Y) - J^2N_J(X,Y)) \notag \\
& = - N_J (X,Y) + N_J (X,Y) \notag \\
& = 0\,.\notag
\end{align}
\end{proof}

The following proposition is needed in Section $2$.

\begin{prop} \label{prop:lieJ}
On an almost complex manifold $(M,J)$, the Lie derivative of $J$ with
respect to a vector field $Y\in\Om^0 (TM)$ is given by
\[
\Ll_Y J = (2J) (\db Y) + \frac{1}{2} J (Y \lrcorner N_J) \in
\Om^{0,1}_J (TM)\,.
\]
\end{prop}
\begin{proof}
We have that
\begin{align}
(\Ll_Y J) (X) & = [Y,JX] - J[Y,X] = J[X,Y] - [JX,Y] \notag \\
& = (2J) \frac{1}{2}\left([X,Y] + J[JX,Y] + 
\frac{1}{2}(N_J(X,Y) - N_J(X,Y)) \right) \notag \\ 
& = (2J)(\db_X Y) + \frac{1}{2} J N_J (Y,X) \,. \notag
\end{align}
\end{proof}

As we will now see, the $\db$-operator on $\Om^{0,1}_J(TM)$ can be 
identified with an appropriate derivative of the Nijenhuis tensor.

\begin{defn} \label{defn:db01TM}
Define the operator
$\db:\Om^{0,1}_J(TM) \rightarrow \Om^{0,2}_J(TM)$, $A \mapsto \db A$,
as the unique linear operator which is given on elements of the form
$A = \al\otimes Z \in \Om^{0,1}_J(M)\otimes\Om^{0} (TM) =
\Om^{0,1}_J(TM)$ by
\[
\db A = \db (\al\otimes Z) = \db\al \otimes Z - \al \wedge \db Z\,.
\]
\end{defn}

Let $\Jj$ denote the space of almost-complex structures on the
manifold $M$. The Nijenhuis tensor can be seen as a map
\[
N:\Jj\to\Om^2 (TM)\,,
\]
with derivative
\[
dN:T\Jj\to\Om^2 (TM)\,.
\]
Given $J\in\Jj$ we have that
\[
T_J \Jj = \left\{A\in\Aut(TM)\,:\ AJ+JA = 0 \right\}
\equiv \Om_J^{0,1}(TM)\,,
\]
which means that
\[
dN_J : \Om_J^{0,1}(TM) \to \Om^2 (TM)\,.
\]
The following lemma follows from a standard calculation.
\begin{lemma} \label{lem:dN}
Given $J\in\Jj$ and $A\in\Om_J^{0,1}(TM)$, we have that
$dN_J(A)\in\Om^2 (TM)$ is given by
\[
dN_J(A)(X,Y) = [AX,JY] + [JX,AY] - J([X,AY]+[AX,Y]) 
- A([X,JY] + [JX,Y])\,,
\]
for any $X,Y\in\Om^0 (TM)$.
\end{lemma}

\begin{prop} \label{prop:dN}
Given $J\in\Jj$ and $A\in\Om_J^{0,1}(TM)$, we have that
$dN_J(A)\in\Om^2 (TM) = 
\Om_J^{2,0}(TM)\oplus\Om_J^{0,2}(TM)\oplus\Om_J^{1,1}(TM)$ can be decomposed
as 
\[
dN_J(A) = (dN_J(A))^{0,2} + (dN_J(A))^{2,0} + (dN_J(A))^{1,1}\,,
\]
where
\begin{align}
(dN_J(A))^{0,2}(X,Y) & = (-2J) (\db A) (X,Y) \,,\notag \\
(dN_J(A))^{2,0}(X,Y) & = \frac{1}{2} J A (N_J(X,Y)) \quad\text{and}
\notag \\ 
(dN_J(A))^{1,1}(X,Y) & = \frac{1}{2} J [N_J(X,AY) + N_J(AX,Y)]\,,\notag
\end{align}
for any $X,Y\in\Om^0 (TM)$.
\end{prop}
\begin{proof}
All the expressions have the right $(p,q)$-type and the $\db$-operator
defined by $J(dN_J - (dN_J)^{2,0} - (dN_J)^{1,1})/2$ has the characterizing
property of Definition~\ref{defn:db01TM}.
\end{proof}
\begin{remark} \label{rem:dN02}
When $J\in\I \subset\Jj$ is an integrable complex
structure, Proposition~\ref{prop:dN} tells us that
\[
dN_J = (dN_J)^{0,2} = (-2J) \db
\]
\end{remark}
\noindent
Let $\Om^{0,2} (TM)$ denote the vector bundle over $\Jj$ whose fiber
over a point $J\in\Jj$ is given by
\[
\Om^{0,2}(TM)|_J = \Om^{0,2}_J (TM)\,.
\]
Since $\Om^{0,2} (TM)$ is a canonical summand of the trivial bundle
$\Om^2(TM) \times \Jj$ over $\Jj$, it carries a natural connection
$\nabla$ defined by projection:
\[
\nabla\cdot = (d\cdot)^{0,2}\,.
\]
The Nijenhuis tensor $N$ can be seen as a natural section of this
vector bundle: 
\[ N : \Jj \to \Om^{0,2}(TM) \,.\]
Proposition~\ref{prop:dN} immediately implies the following
generalization of Remark~\ref{rem:dN02}.
\begin{cor} \label{cor:dN02}
For any $J\in\Jj$ we have that
\[
\nabla N_J = (dN_J)^{0,2} = (-2J) \db \,.
\]
\end{cor}

\section{A commutation relation for K\"ahler manifolds}
\label{app:C}
\noindent
Our goal here is to prove that on a K\"ahler manifold $(M,J,\om)$, 
with Riemannian metric given by $g(\cdot,\cdot) = \om(\cdot,J\cdot)$, 
the diagram
\[ 
\xymatrix{  X\in \Omega^0(TM) \ar[r]^{\bar{\partial}} \ar[d] 
          & \Omega^{0,1}(TM)\ni\Ga \ar[d] \\
            \al_{\Tiny X}\in\Omega^{0,1}(M) \ar[r]^{\bar{\partial}} 
          & \Omega^{0,2}(M)\ni\al_{\Tiny\Ga} } 
\]
commutes, where
\[
\al_{\Tiny X}(Y) = g(X,Y) - i\om(X,Y)
\]
and
\[
\al_{\Tiny\Ga} (Y,Z) = [g(\Ga(Y),Z) - g(\Ga(Z),Y)]
- i [\om(\Ga(Y),Z) - \om(\Ga(Z),Y)]\,.
\]
In other words, we need to show that
\begin{equation} \label{eq:claim}
\db (\al_{\Tiny X}) - \al_{\Tiny \db X} = 0\,,\ \forall X\in
\Omega^0(TM) \,.
\end{equation}

\begin{lemma} \label{lem:claim}
The map
\begin{align}
\Om^0(TM) & \to \Om_J^{0,2} (M) \notag \\
X & \mapsto \db (\al_{\Tiny X}) - \al_{\Tiny \db X} \notag
\end{align}
is a tensor.
\end{lemma}
\begin{proof}
As usual, it suffices to show that the given map is
$\Om^0(M)$-linear. For any function $f\in\Om^0(M)$ and vector field
$X\in\Om^0(TM)$, we have that
\[
\db(\al_{\Tiny fX}) = \db (f \al_{\Tiny X}) = (\db f) \wedge
\al_{\Tiny X} + f \db \al_{\Tiny X}
\]
while
\[
\db(fX) = (\db f) \otimes X + f \db X 
\quad\Rightarrow\quad
\al_{\Tiny \db(f X)} = \al_{\Tiny (\db f) \otimes X} + f
\al_{\Tiny\db X}\,. 
\]
Since
\begin{align}
\al_{\Tiny (\db f) \otimes X} (Y,Z) 
& = [g((\db f)(Y)X,Z) - g((\db f)(Z)X,Y)]
\notag \\
& \quad - i [\om ((\db f)(Y)X,Z) - \om ((\db f)(Z)X,Y)] 
\notag \\
& = (\db f)(Y) [g(X,Z) - i \om(X,Z)] 
\notag \\
& \quad - (\db f)(Z) [g(X,Y) - i \om(X,Y)] 
\notag \\
& = (\db f)(Y) \al_{\Tiny X}(Z) -  (\db f)(Z) \al_{\Tiny X}(Y)
\notag \\
& = \left( \db f \wedge \al_{\Tiny X}\right) (Y,Z) \,,
\notag
\end{align}
we conclude that
\begin{align}
\db(\al_{\Tiny fX}) - \al_{\Tiny \db(fX)} 
& = \left( \db f \wedge \al_{\Tiny X} + f
  \db\al_{\Tiny X}\right) -  \left( \db f \wedge \al_{\Tiny X} +
  f \al_{\Tiny \db X} \right) \notag \\
& = f \left( \db (\al_{\Tiny X}) - \al_{\Tiny \db X} \right)\,, \notag 
\end{align}
which proves the lemma.
\end{proof}
\noindent
Lemma~\ref{lem:claim} implies that it suffices to
prove~(\ref{eq:claim}) at an arbitrary point. Since a K\"ahler
manifold behaves up to first order at a point as flat $\C^n$, 
where~(\ref{eq:claim}) clearly holds, we can conclude
that~(\ref{eq:claim}) is true on any K\"ahler manifold.

\section{The $A_\infty$-action of $G_\lambda$ on the tubular neighborhoods}
\label{quasi}

In this section we fill in some details in the proof of Theorem \ref{maintopological}.

We will need to use the notion of homotopy (co)limits. We have already used two examples of homotopy colimits: the homotopy pushout \eqref{hopdefn} and the homotopy orbits of a group action \eqref{hoorbdefn}. The homotopy colimit of a sequence of maps is also familiar: it is the infinite mapping telescope. A friendly and elegant reference for homotopy limits and colimits is \cite{HV} (see also 
\cite{BK} and \cite{V}). 

As usual, we regard posets as categories with exactly one arrow between two objects $a$ and $b$ when 
$a\leq b$.

\medskip

\noindent

If $X$ is a space, we denote by $\cat{P}_X$ the poset of 
subspaces of $X$ ordered by \emph{reverse} inclusion.

Let $G$ be a topological group. We denote by $\K_G$ the partially ordered set of compact subspaces of $G$ ordered by inclusion. There is a canonical map 
\[  A_G = \hocolim_{K \in \K_G} K \xrightarrow{\phi} G \]
which is a weak homotopy equivalence. 

There is a strictly associative multiplication on $A_G$ induced by the functor 
\begin{eqnarray*}
  \K_G \times \K_G & \to  & \K_G \\
  (K,L) & \mapsto & KL 
\end{eqnarray*}
and the natural isomorphism \cite[Proposition 3.1(4)]{HV} 
\[ A_G \times A_G = \hocolim_{(K,L) \in \K_G \times \K_G } K\times L. \]
Note that $\phi$ is a strictly multiplicative weak equivalence.

Consider the set
\[ \KS_G = \{ (K_n,\ldots,K_1) \mid n\geq 0, \quad K_i \subset G \text{ compact and } 
\neq \emptyset\} \]
(when $n=0$ we mean the empty word) with the partial order defined by 
\[ (K_n,\ldots,K_1) \leq (H_m,\ldots,H_1) \text{ if } n \leq m \text{ and } K_i \subset H_i. \]
Given sequences $S,T \in \KS_G$ we write $S\ast T$ for their concatenation.

\begin{defn}
Let $G$ be a topological group, $X$ a $G$-space and $U\subset X$ a subspace (not necessarily
$G$-invariant). A \emph{near action} of $G$ on $U$ consists of a functor
\begin{eqnarray*}
\KS_G  & \to & \cat{P}_U \\
(K_n,\ldots,K_1) & \mapsto & U_{(K_n,\ldots,K_1)} 
\end{eqnarray*}
such that
\begin{enumerate}[(i)]
\item $U_\emptyset = U$ (where $\emptyset$ denotes the empty sequence),
\item For each $S \in \cat K$, the inclusion $U_S \to U$ is a weak equivalence,
\item Given $K \in \K_G$ and $T \in \KS_G$, 
the restriction of the $G$-action
\[ K \times U_{(K)*T} \to X \]
has image contained in $U_T$.
\end{enumerate}
A \emph{near $G$-equivariant map} is a natural transformation of functors commuting 
strictly with the action of the compact subsets.
\end{defn}

\begin{lemma}
\label{repaction}
A near action of $G$ on $U$ induces a canonical action of $A_G$ on 
$T(U) = \holim_{S \in \KS_G} U_S$.
\end{lemma}
\begin{proof}
There is an obvious action $K\times T(U) \to T(U)$ for each compact set and this extends
canonically to the required action.
\end{proof}

Note that there is a canonical homotopy equivalence 
\[ T(U) \xrightarrow{\pi} U \]
induced by the inclusion of the empty sequence in $\KS_G$.

\begin{remark}
An $A_\infty$-action (see \cite{St}) of a topological monoid $G$ on a space $X$ is a map $G\times X \to X$ which is "associative up to all higher homotopies" with respect to the multiplication on $G$. One way of giving such an $A_\infty$-action is to give an actual action of a topological monoid $G'$ on a space $X'$ together with weak equivalences $G \to G'$ and $X \to X'$. One can then perform homotopy meaningful constructions (such as homotopy orbits) replacing $G,X$ with $G',X'$. 

Setting $G'=A_G$ and $U'=T(U)$, the previous Lemma shows that a near action of $G$ on $U$ gives rise to an $A_\infty$-action of $G$ on $U$.
\end{remark}

\begin{defn}
\label{defhorbit}
The \emph{homotopy orbit space} $\overline{U_{hG}}$ of a near action of $G$ on $U$ is the realization of the semi-simplicial space 
\[  n \mapsto  A_G^n\times T(U) \]
\end{defn}

A $G$-space $U$ has a trivial near $G$-action where $U_S=U$ for all $S \in \KS_G$ and 
unless we specify otherwise we always give $G$-spaces this near $G$-action.
\begin{lemma}
\label{equivnear}
If $U$ is $G$-invariant there is a weak equivalence 
\[  \overline{U_{hG}} \to U_{hG}. \]
\end{lemma}
\begin{proof}
For each compact set there is a canonical homotopy making the diagram 
\[ \xymatrix{ K \times T(U) \ar[r] \ar[d] &  K \times U \ar[d] \\ T(U) \ar[r] & U }\]
and this yields a homotopy coherent weak equivalence between the two semi-simplicial spaces
in question. The result follows by taking realization.
\end{proof}

\noindent
We can now elaborate on the proof of Theorem \ref{maintopological}. With reference to the 
notation in the statement of that theorem, in the remainder of this
section we will write $K$ for the isometry group $K(m)$, $W$ for the representation
of $K$ on the normal to the corresponding stratum (described in Theorem \ref{representation}),
$V$ for the stratum $V_{\ell+1}$, $NV$ for its tubular neighborhood, $X$ for $\cat I_{\lambda+1}$
and $G$ for $G_{\lambda+1}$.

We'll fix a $K$-invariant metric on $W$ and given a continuous function 
$\epsilon \colon G \to \R_+$ write 
\[ G  \times_K (W \setminus 0)_\epsilon = 
\{ g \cdot w \mid g \in G, 0< |w| < \epsilon(g) \} \]

\begin{prop}
\label{identtube}
There is a continuous function $\epsilon: G \to \R_+$ such that 
\[ \xymatrix{ G \times_K (W \setminus 0)_\epsilon \ar@{^{(}->}[r]^{\quad \psi} \ar[d]  & \ar[d]^\pi NV\setminus V\\
 G/K \ar@{^{(}->}[r] & V } \]
commutes in the (weak) homotopy category. Moreover $\psi$ is a weak equivalence.
\end{prop}
\begin{proof}
The function $\epsilon$ exists by continuity of the action.
The slice theorem for the action of the symplectomorphism group $G$ on the space of compatible almost complex structures\footnote{The construction of the slice for the action of the diffeomorphism group on the space of metrics in \cite{Eb} works in this case.} together with the uniqueness of tubular neighborhoods give, for each right $K$-invariant compact subset $L \subset G$,
a homeomorphism $\psi_L$, homotopic to the inclusion, such that the diagram
\[ \xymatrix{
L \times_K (V\setminus 0) \ar[rd] \ar@{^{(}->}[rr]^{\psi_L} & & \pi^{-1}(L/K) \ar[ld]^\pi \\ & L/K & }
\]
commutes. The uniqueness of tubular neighborhoods implies that if $L \subset L'$, ${\psi_{L'}}_{|(L\times_K (V\setminus 0))}$ is homotopic over $L/K$ to $\psi_L$. 
The result follows.
\end{proof}

\begin{remark}
The crucial point in the proof of the previous proposition is the existence of a slice for the 
action. One can apply the arguments in this section whenever this is the case. 
\end{remark}

The subspace 
\[ G \times_K (W\setminus 0)_\epsilon \subset NV\setminus V \]
can be endowed with a near $G$-action, by choosing for each sequence $S \in \KS_G$
a continuous function $\epsilon_S \colon G \to \R_+$ with $\epsilon_\emptyset= \epsilon$ 
in such a way that for each compact subset $L \subset G$,
\[ L \cdot (G \times_K (W\setminus 0)_{\epsilon_{(L)*S}})  \subset 
G \times_K (W\setminus 0)_{\epsilon_{S}}. \]

Giving a $G$-space $U$ the trivial near $G$-action we have a 
pushout diagram of near $G$-spaces and near $G$-equivariant maps
\[  G/K \longleftarrow G \times_K (W \setminus 0)_\epsilon \longrightarrow (X\setminus V). \]
Writing $P$ for the homotopy pushout, there is an obvious near $G$-action on $P$ together
with a near $G$-equivariant map 
\[ P \to X \]
which is clearly a weak equivalence. Since the canonical map 
\[ \hocolim \left( T(G/K) \leftarrow T(G \times_K (W \setminus 0)_\epsilon) \rightarrow T(X\setminus V) \right) \to T(P) \]
is a weak equivalence, applying $A_G$ homotopy orbits (in the sense of 
Definition \ref{defhorbit}) and Lemma \ref{equivnear} we get
\begin{equation}
\label{nearpushout}
\hocolim \biggl( BK \leftarrow (W\setminus 0)_{hK} \rightarrow (X\setminus V)_{hG} \biggr)
\simeq X_{hG} 
\end{equation}
as required.

\end{document}